\documentclass{amsart}

\usepackage{amsmath,amsthm,amssymb,amsfonts}
\usepackage{bbm}
\usepackage[mathscr]{eucal}
\usepackage[english]{babel}

\allowdisplaybreaks[1]

\newtheorem{thm}{Theorem}[section]

\newtheorem{lem}[thm]{Lemma}
\newtheorem{prop}[thm]{Proposition}
\newtheorem{rem}[thm]{Remark}

\newtheorem{defn}[thm]{Definition}

\theoremstyle{remark}

\numberwithin{equation}{section}

\newcommand{\thmref}[1]{Theorem~\ref{#1}}
\newcommand{\lemref}[1]{Lemma \ref{#1}}
\newcommand{\propref}[1]{Proposition \ref{#1}}

\def\CC{{\mathbb C}}
\def\NN{{\mathbb N}}
\def\RR{{\mathbb R}}
\def\ZZ{{\mathbb Z}}

\def\cQ{\mathcal{Q}}

\def\cS{\mathcal{S}}

\def\cX{\mathcal{X}}

\def\fD{{\mathfrak D}}

\def\ff{{\mathfrak f}}
\def\fg{{\mathfrak g}}
\def\qq{{\mathfrak q}}

\def\sD{\mathscr{D}}

\def\sP{\mathscr{P}}

\def\eps{{\varepsilon}}
\def\QQ{{\cQ}}
\def\ONE{{\mathbbm 1}}
\def\II{\mathcal{I}}
\def\Q{Q}
\def\Om{\Omega}
\def\FF{\Phi}
\def\aa{\beta}
\def\cc{a}
\def\BB{{\dot{B}}}

\def\BMO{\operatorname{BMO}}
\def\VMO{\operatorname{VMO}}
\def\loc{{\operatorname{loc}}}
\def\avg{{\operatorname{Avg}}}
\def\supp{\operatorname{supp}}

\begin{document}

\title[Nonlinear spline approximation in $\BMO(\RR)$]
{Nonlinear spline approximation in $\BMO(\RR)$}

\author[K. G. Ivanov]{Kamen G. Ivanov}
\address{Institute of Mathematics and Informatics\\
Bulgarian Academy of Sciences\\
Sofia, Bulgaria}
\email{kamen@math.bas.bg}

\author[P. Petrushev]{Pencho Petrushev}
\address{Department of Mathematics\\
University of South Carolina\\
Columbia, SC}
\email{pencho@math.sc.edu}

\thanks{The first author has been supported by Grant DN 02/14
of the Fund for Scientific Research of the Bulgarian Ministry of Education and Science.
The second author has been supported by NSF Grant DMS-1714369.}

\begin{abstract}
We study nonlinear approximation in $\BMO$ from splines
generated by a hierarchy of B-splines over regular multilevel nested partitions of $\RR$.
Companion Jackson and Bernstein estimates are established that allow to completely characterize
the associated approximation spaces.
\end{abstract}

\maketitle

\begin{sloppypar}

\section{Introduction}\label{sec:introduction}

The space of functions of bounded mean oscillation ($\BMO$)
is a natural replacement of $L^\infty$ in many problems in Analysis and, in particular, in Approximation theory.
Here we consider nonlinear approximation from regular splines generated by a nested hierarchy of B-spines in $BMO$ on $\RR$.

As usual the space BMO is defined as the set of all functions $f\in L^1_{\loc}(\RR)$ such that
\begin{equation*}
\|f\|_{\BMO}:=\sup_{I}\frac{1}{|I|}\int_I|f(x)-\avg_I f| dx <\infty,
\quad \avg_I f:=\frac{1}{|I|}\int_I f(x)dx,
\end{equation*}
where the $\sup$ is over all intervals $I\subset \RR$.
As is well known the set $C_0(\RR)$ of all continuous functions with compact support is not dense in $\BMO$.
Therefore, it is natural to approximate in the $\BMO$-norm functions that belong to the space $\VMO$,
the closure of the space $C_0(\RR)$ in the $\BMO$-norm, see \cite{CW}.

In this article we consider nonlinear $n$-term approximation from $B$-splines
generated by multilevel nested partitions of $\RR$.
More specifically, let $\{\II_m\}_{m\in\ZZ}$ be a sequence of families of intervals
such that each level $\II_m$ is a partition of $\RR$ into compact intervals with disjoint interiors
and a refinement of the previous level $\II_{m-1}$.
We consider regular partitions of this sort, i.e. we require that the intervals
from each level $\II_m$ be of comparable length.
Define $\II:=\cup_{m\in\ZZ}\II_m$.
Each such multilevel partition $\II$ generates a ladder of spline spaces
$\cdots\subset S_{-1}^k \subset S_0^k\subset S_1^k \subset \cdots$
of degree $k-1$, where $S_m^k$ is spanned by B-splines $\{\varphi_\Q\}$.
We denote by $\QQ_m$ the supports of the $m$th level B-splines $\varphi_\Q$
and set $\QQ:=\cup_{m\in\ZZ} \QQ_m$.

We are interested in approximating from the nonlinear set $\Sigma_n$ of all splines of the form
\begin{equation*}
g=\sum_{\Q\in \Lambda_n} c_\Q\varphi_\Q,
\quad \Lambda_n\subset \QQ, \; \# \Lambda_n \le n.
\end{equation*}
Here the index set $\Lambda _n$ is allowed to vary with $g$.
The approximation error is defined by
\begin{equation*}
\sigma_n(f)_{\BMO}: = \inf_{g\in \Sigma_n} \|f-g\|_{\BMO}.
\end{equation*}

The Besov spaces $\BB_\tau^{\alpha, k}$ with $\alpha>0$, $\tau:= 1/\alpha$, and $k\in \NN$,
play a crucial role here.
In the case when $\tau \ge 1$ the space $\BB_\tau^{\alpha, k}$
consists of all functions $f\in L^\tau_\loc(\RR)$
such that $\Delta_h^kf\in L^\tau(\RR)$, $\forall h\in\RR$, and
\begin{equation*}
\|f\|_{\BB_\tau^{\alpha, k}}\!:=\!\Big(\int_0^\infty\![t^{-\alpha}\omega_k(f, t)_\tau]^\tau\frac{dt}{t}\Big)^{1/\tau}\!<\infty,
~~~
\omega_k(f, t)_\tau\!:=\!\sup_{|h|\le t} \|\Delta_h^kf(\cdot)\|_{L^\tau(\RR)}.
\end{equation*}
In the case when $\tau<1$ the above definition the Besov space $\BB_\tau^{\alpha, k}$
is not quite satisfactory to us.
We believe that in principle the use of $\omega_k(f, t)_\tau$ with $\tau<1$, as  well as  polynomial approximation in $L^\tau$, $\tau<1$,
should be avoided when possible.
For this reason we define the Besov space $\BB_\tau^{\alpha, k}$ when $\tau<1$
via local polynomial approximation in $L^q$, $q\ge 1$,
(see Definition~\ref{def:B-2}).

Clearly, $\|f+P\|_{\BB_\tau^{\alpha, k}}=\|f\|_{\BB_\tau^{\alpha, k}}$
for all $P\in\Pi_k$, the set of all algebraic polynomials of degree $k-1$.
Hence, $\BB_\tau^{\alpha, k}$ consists of equivalence classes modulo $\Pi_k$.
As will be shown (Theorem~\ref{thm:rep-B}) for any function $f\in \BB_\tau^{\alpha, k}$
there exists a polynomial $P\in\Pi_k$ such that $f-P\in \VMO$.
We shall be assuming that each $f\in \BB_\tau^{\alpha, k}$ is the {\em canonical representative}
$f-P\in \VMO$ of the equivalence class modulo $\Pi_k$ generated by $f$.

Our primary goal in the article is to prove the following Jackson and Bernstein estimates (Theorems~\ref{thm:Jackson} and \ref{thm:bernstein}):
Let $\alpha>0$, $\tau:=1/\alpha$, $k\ge 2$.
If $f\in \BB_\tau^{\alpha, k}$,  then $f\in \VMO$ and
\begin{equation*}
\sigma_n(f)_{\BMO} \le cn^{-\alpha}\|f\|_{\BB_\tau^{\alpha, k}}, \quad n\ge 1,
\end{equation*}
and for any $g\in \Sigma_n$
\begin{equation*}
\|g\|_{\BB_\tau^{\alpha, k}} \le cn^\alpha\|g\|_{\BMO}.
\end{equation*}
As is well known these two estimates imply a complete characterzation
of the approximation spaces associated to this sort of approximation
(see Theorem~\ref{thm:app-sp} below).

To achieve our objectives we first develop in sufficient detail the theory of the Besov spaces $\BB_\tau^{\alpha, k}$.
In particular we derive representations of these spaces in terms of
B-splines, local polynomial approximation, and spline quasi-interpolants.

We also compare the nonlinear spline approximation in $\BMO$
with the spline approximation in the uniform norm and in $L^p$, $1\le p<\infty$.

There is a considerable difference between nonlinear spline approximation in $\BMO$
in dimension $d=1$ and in dimension $d>1$.
A detailed discussion and clarification of this phenomenon will be given
in a followup article.

Observe that nonlinear approximation in $\BMO(\RR^d)$ from dyadic piecewise polynomial functions
has been studied by I. Irodova, see \cite{Irodova} and the references therein.
In this case dyadic Besov spaces are used for characterization of the rates of approximation.

There is a close relationship between nonlinear approximation from splines and wavelets in $\BMO$.
We develop the nonlinear $n$-term wavelet approximation in $\BMO$ in \cite{IP3}.

\smallskip
\noindent
{\bf Outline.}
This article is organized as follows.
In Section~\ref{sec:prelim} we introduce our setting and collect all facts we need
concerning splines, local polynomial approximation, spline quasi-interpolants.
The Besov spaces involved in nonlinear spline approximation in $\BMO$ are studied in Section~\ref{sec:Besov-spaces}.
In Section~\ref{sec:spline-approx} we state our main results on spline approximation in $\BMO$
and we compare them in Section~\ref{sec:comparison}
with the existing results on spline approximation in the uniform norm and in $L^p$, $1\le p<\infty$.
Sections~\ref{sec:jackson} and \ref{sec:bernstein} contain the proofs of our Jackson and Bernstein estimates.
Section~\ref{sec:appendix} is an appendix, where we place the proofs of some statements from previous sections.

\smallskip
\noindent
{\bf Notation.}
We shall use the notation $\|\cdot\|_p:=\|\cdot\|_{L^p(\RR)}$.
$C_0(\RR)$ will stand for the set of all continuous and compactly supported functions on $\RR$.
Given a measurable set $A\subset \RR$ we denote by $|A|$ its Lebesgue measure
and by $\ONE_A$ its characteristic function.
For an interval $J$ and $\lambda>0$ we denote by $\lambda J$ the interval of length $\lambda|J|$ which is co-centric with $J$.
As usual $\ZZ$ will denote the set of all integers,
$\NN$ will be the set of all positive integers, and $\NN_0:=\NN\cup\{0\}$.
Also, $\Pi_k$ will stand for the set of all univariate algebraic polynomials of degree $k-1$.
Unless specified otherwise, all functions are complex-valued.
Positive constant will be denoted by $c, c', c_1, \dots$ and they may vary at every occurrence;
$a\sim b$ will stand for $c_1\le a/b\le c_2$.

\section{Preliminaries}\label{sec:prelim}

In this section we collect all facts regarding the $\BMO$ space, B-splines,
local approximation, quasi-interpolants, and other results that will be needed in the sequel.

\subsection{$\BMO$ and $\VMO$ spaces}\label{subsec:BMO}

As usual the space BMO on $\RR$ is defined as the set of all locally integrable functions on $\RR$
such that
\begin{equation}\label{def-BMO-J}
\|f\|_{\BMO}:=\sup_{J}\frac{1}{|J|}\int_J|f(x)-\avg_{J} f| dx <\infty,
\quad \avg_J f:=\frac{1}{|J|}\int_J f(x)dx,
\end{equation}
where the $\sup$ is over all compact intervals $J$ and $|J|$ is the length of $J$.
Note that $\|\cdot\|_{\BMO}$ is not a norm
because $\|g\|_{\BMO}=0$ if $g= {\rm constant}$ ($g\in\Pi_1$).
For this reason we identify each $f\in\BMO$ with $f+a$, $a =$ constant,
and view $\BMO$ as a subset of $L^1_{\loc}(\RR)/\Pi_1$.
Then $\|\cdot\|_{\BMO}$ is a norm on $\BMO$.

From the well known John-Nirenberg inequality \cite{JN} it follows that
any function $f\in \BMO$ is in $L^p_{\loc}(\RR)$, $1<p<\infty$, and
\begin{equation}\label{equiv-norms-BMO}
\sup_{J}\Big(\frac{1}{|J|}\int_J|f(x)-\avg_J f|^p dx\Big)^{1/p}
\sim \sup_{J}\frac{1}{|J|}\int_J|f(x)-\avg_J f| dx
=\|f\|_{\BMO}
\end{equation}
with constants of equivalence depending only on $p$.

As is well known $\BMO$ is \emph{not a separable space}.
However, the space $\VMO$, defined as the closure of $C_0(\RR)$ in the $\BMO$ norm
(see Section~4 in \cite{CW}),
is a \emph{separable} Banach space.
As in $\BMO$ the elements of $\VMO$ are classes of equivalence modulo constants.
We shall consider spline approximation of functions in $\VMO$.

\subsection{Nested partitions of $\RR$}\label{subsec:nested_partitions}

We say that
\begin{equation*}
\II = \cup_{m\in\ZZ} \II_m
\end{equation*}
is a \emph{regular multilevel partition} of $\RR$ into subintervals
with levels $\{\II_m\}$ if the following three conditions are obeyed:

(a) Each level $\II_m$ is a partition of $\RR$, i.e. $\RR=\cup_{I\in \II_m} I$,
and $\II_m$ consists of compact intervals with disjoint interiors.

(b) The levels $\{\II_m\}$ of $\II$ are nested, i.e. $\II_{m+1}$ is a refinement of $\II_m$,
and each $I\in \II_m$ has at least two and at most $M_0$ children in $\II_{m+1}$,
where $M_0\ge 2$ is a constant independent of $m$.

(c) There exists a constant $\lambda\ge 1$ such that
\begin{equation}\label{cond-rho}
|I'|\le \lambda |I''|,
\quad\forall I', I''\in \II_m, \;\;\forall m\in\ZZ.
\end{equation}

\smallskip

The dyadic intervals $\sD= \cup_{m\in\ZZ} \sD_m$, $\sD_m=\{[2^{-m}k,2^{-m}(k+1)] : k\in\ZZ\}$,
are an example of a regular multilevel partition of $\RR$ with $\lambda=1$ and $M_0=2$.

\smallskip
\noindent
{\bf Other scenarios.}
Another version of the above setting is when the multilevel partition $\II$ of $\RR$ is of the form
$\II = \cup_{m=0}^\infty \II_m$,
where the levels $\II_m$ satisfy conditions (a)--(c) from above.

Yet a third version of this setting is
when we have a regular multilevel
partition $\II = \cup_{m\ge 0} \II_m$ of a fixed compact interval $J$
obeying conditions similar to (a)--(c) above and $J$ is the only element of $\II_0$.

Our theory can readily be adjusted to each of these settings.
We shall stick to the first setting on $\RR$ from above.

\smallskip

A couple of remarks are in order.

(1) There is a natural tree structure in $\II$ induced by the inclusion relation.

\smallskip

(2) Conditions (a)--(c) imply that there exist constants $0<r<\rho<1$
such that if $I\in\II_m$, $I'\in\II_{m+1}$, and $I'\subset I$, then
\begin{equation}\label{r-rho}
r\le |I'|/|I|\le \rho.
\end{equation}
In fact, inequality \eqref{r-rho} is derived with bounds $1/(M_0\lambda-\lambda+1)\le r\le 1/2$, $\rho\le\lambda/(\lambda+1)$.
Note that the upper bound for $\rho$ and the lower bound for $r$ cannot be simultaneously achieved provided $M_0\ge 3$.

We denote by $x_{m,j}, j\in \ZZ$, the knots of $\II_m$, i.e. $\II_m\!=\!\{[x_{m, j}, x_{m,j+1}]: j\!\in\!\ZZ\}$,
where
$$\cdots < x_{m, -1} < x_{m, 0}< x_{m, 1} < \cdots.$$
With this notation the nested condition (b) can be described as
\begin{equation*}
\{x_{m,j} : j\in \ZZ\}\subset \{x_{m+1,j} : j\in \ZZ\},\quad m\in\ZZ.
\end{equation*}
Observe that inequality \eqref{r-rho} implies that every regular multilevel partition $\II$ may have either zero or one common knot for all levels $\II_m$.
Dyadic intervals $\sD$ are example with one common knot -- the origin.

For a fixed integer $k\ge 2$ we denote by $\QQ_m$ the collection of all unions of $k$ consecutive intervals from $\II_m$, i.e.
\begin{equation}\label{def-Q}
\QQ_m=\{\Q=[x_{m,j}, x_{m, j+k}] : j\in\ZZ\},\quad m\in\ZZ.
\end{equation}
Further, set $\QQ(\II)=\QQ:=\cup_{m\in\ZZ} \QQ_m$.
The intervals from $\QQ$ are the supports of the $B$-splines of degree $k-1$ defined below in \S\ref{subsec:splines}.
We shall use both $\II$ and $\QQ$ as index sets of the objects discussed in this article.

To every $I\in\II$ we associate the set
\begin{equation}\label{def-Omega-I}
\Omega_I:=\cup \{\Q: \Q\in \QQ_m, \Q\supset I\}, \quad I\in \II_m.
\end{equation}
In other words, if $I=[x_{m,j}, x_{m, j+1}]$ then $\Omega_I=[x_{m,j+1-k}, x_{m, j+k}]$.
In view of condition (c) of the regular multilevel partition $\II$ we always have
\begin{equation*}
|I|\sim|\Q|\sim|\Omega_I|,\quad \forall I\in\II_m, \Q\in\QQ_m, m\in\ZZ.
\end{equation*}

Given a regular multilevel partition $\II$, we define \emph{the level of a compact interval $J$} as the largest $\nu\in\ZZ$
such that $\II_\nu$ has no more than one knot in the interior of $J$.
Observing that $\II_{\nu+1}$ has at least two knots in the interior of $J$ we infer from \eqref{r-rho} that $|J|\sim |I|$ for every $I\in\II_\nu$.
Also for every $m\le \nu$ we have two adjacent intervals $I_1$ and $I_2$ in $\II_m$ such that $J\subset I_1\cup I_2$ and $|I_1|\sim|I_1|\ge c|J|$.

\subsection{Splines over nested partitions of $\RR$}\label{subsec:splines}

\smallskip
\noindent
{\bf Piecewise polynomials and splines.}
Assuming that $\II$ is a regular multilevel partition of $\RR$,
we denote by $\tilde\cS_m^k := \tilde\cS^k(\II_m)$, $m \in \ZZ$, $k \ge 1$,
the set of all piecewise polynomial (possibly discontinuous) functions over $\II_m$ of degree $k-1$,
i.e.
\begin{equation*}
S \in \tilde\cS_m^k
\quad\hbox{if}\quad
S = \sum_{I \in \II_m} \ONE_I \cdot P_I,
\end{equation*}
where $\ONE_I$ is the characteristic function of $I$ and $P_I \in \Pi_k$.
$S$ is assumed to be right-continuous at the knots of $\II_m$.
Then the space of the $m$th level splines is defined by
\begin{equation}\label{def-S}
\cS_m^k := \tilde\cS_m^k\cap C^{k-2}, \quad k\ge 2.
\end{equation}

\smallskip
\noindent
{\bf B-splines.}
Given $\Q:=[x_{m, j}, x_{m, j+k}]\in\QQ_m$, $m\in\ZZ$, (see \eqref{def-Q}) we denote by $\varphi_\Q$
the $B$-spline of degree $k-1$ with knots $x_{m, j}, \dots, x_{m, j+k}$; these are $k+1$ consecutive knots of $\II_m$.
For the precise definition of $\varphi_\Q$, see e.g. \cite[Chapter~5, (2.7)]{DL}.
Note that: (i) $\Q$ is the support of $\varphi_\Q$ and
(ii) $\|\varphi_\Q\|_\infty\sim 1$ with constants of equivalence depending only on $k$ and $\lambda$.
We denote by $V_\Q:=\{x_{m, j}, \dots, x_{m, j+k}\}$ the knots of $\varphi_\Q$.

It will be convenient to index the B-splines of degree $k-1$ by their supports.
Thus, given a regular multilevel partition $\II$,  the collection of all B-splines of degree $k-1$ is
\begin{equation*}
\Phi=\Phi(\II):= \{\varphi_\Q: \Q\in\QQ(\II)\}.
\end{equation*}

As is well known $\{\varphi_\Q: \Q\in\QQ_m\}$ is a basis for $\cS_m^k$.
Each $S\in \cS_m^k$ has a unique representation
\begin{equation*}
S=\sum_{\Q\in\QQ_m} \aa_\Q(S)\varphi_\Q.
\end{equation*}
According to de Boor -- Fix theorem (see e.g. \cite[Chapter~5, Theorem 3.2]{DL}) the coefficient $\aa_\Q(S)$ is a linear functional given by:
\begin{equation}\label{a-Q}
\aa_\Q(S):= \sum_{\nu=0}^{k-1} (-1)^\nu \varpi_\Q^{(k-\nu-1)}(\xi_\Q)S^{(\nu)}(\xi_\Q),
\end{equation}
\begin{equation*}
\varpi_\Q(x):=\frac{1}{(k-1)!} \prod_{\nu=j+1}^{j+k-1} (x-x_{m, \nu}),\quad \Q=[x_{m, j}, x_{m, j+k}],
\end{equation*}
where $\xi_\Q$ is an arbitrary point from $(x_{m, j}, x_{m,j+k})$.
The value of $\aa_\Q(S)$ in \eqref{a-Q} is independent of the choice of $\xi_\Q$.

From condition (c) on the multilevel partition $\II$ it readily follows that
\begin{equation*}
|\aa_\Q(S)|\le c\|S\|_{L^\infty(\Q)},
\quad \Q\in\QQ_m.
\end{equation*}
This implies that (see e.g. \cite[Lemma~2.3]{DP} in the case $\RR^2$)
if $S=\sum_{\Q\in\QQ_m}b_\Q\varphi_\Q$, where $\{b_\Q\}$ is an arbitrary sequence of complex numbers,
then 
\begin{equation*}
\|S\|_p \sim \Big(\sum_{\Q\in\QQ_m} \|b_\Q\varphi_\Q\|_p^p\Big)^{1/p},
\quad 0<p\le \infty.
\end{equation*}
Moreover, for any $0<p, \tau\le \infty$
\begin{equation}\label{norm-spline}
\Big(\sum_{I\in\II_m} \big[\|S\|_{L^p(I)}\big]^\tau\Big)^{1/\tau}
\sim \Big(\sum_{\Q\in\QQ_m} \big[\|b_\Q\varphi_\Q\|_p\big]^\tau\Big)^{1/\tau}
\end{equation}
with the usual modification when $\tau=\infty$.
Note that $\aa_\Q(S)=b_\Q$.

\subsection{Local polynomial approximation}\label{local-approx}

Recall first some simple properties of polynomials that will be frequently used.

\begin{lem}\label{lem:poly-norms}
Let $P \in \Pi_k$, $k \ge 1$, and $0< p,q \le \infty$.

$(a)$ For any compact interval $J$
\begin{equation}\label{equiv-norms}
|J|^{-1/p}\|P\|_{L^p(J)} \sim |J|^{-1/q} \|P\|_{L^q(J)}
\end{equation}
with constants of equivalence depending only on $k$ and $\min\{p,q\}$.

$(b)$
Let $J^\prime \subset J$ be two intervals such that
$|J| \le (1+\delta)|J^\prime|$, $\delta>0$. Then
\begin{equation*}
\|P\|_{L^p(J)} \le c \|P\|_{L^p(J^\prime)}
\end{equation*}
with $c= c(p,k,\delta)$.

$(c)$
For any compact interval $J$
\begin{equation*}
\|P'\|_{L^p(J)} \le c|J|^{-1} \|P\|_{L^p(J)}
\end{equation*}
with $c=c(p, k)$.
\end{lem}

Applied to B-splines \lemref{lem:poly-norms}~(a), (c) imply
\begin{equation*}
|Q|^{-1/q} \|\varphi_\Q\|_{q}\sim \|\varphi_\Q\|_\infty\sim 1,\quad \|\varphi_\Q'\|_p\sim |Q|^{-1} \|\varphi_\Q\|_p.
\end{equation*}
For a function $f \in L^p(J)$, defined on a compact interval $J \subset\RR$,
$1\le p\le\infty$, and $k \ge 1$, we define
\begin{equation}\label{Ek-J}
E_k(f, J)_p := \inf_{P \in \Pi_k} \|f-P\|_{L^p(J)}.
\end{equation}
Also, we denote by $\omega_k(f, J)_p$ the $k$-th modulus of
smoothness of $f$ on $J$:
\begin{equation}\label{omega-J}
\omega_k(f,J)_p := \sup_{h \in \RR}\|\Delta_h^k (f,\cdot,J)\|_{L^p(J)},
\end{equation}
where
\begin{equation}\label{delta-f}
\Delta_h^k (f, x, J):=\left\{
\begin{array}{l l}
\sum_{j=0}^k (-1)^{k+j} \binom{k}{j} f(x+jh),&\mbox{if}~[x, x+kh]\subset J;\\
0,&\mbox{otherwise}.
\end{array}
\right.
\end{equation}

\begin{lem}[Whitney]\label{lem:whitney}
If $f \in L^p(J)$,
$1\le p \le \infty$, and $k \ge 1$, then
\begin{equation*}
E_k(f, J)_p \le c\omega_k(f, J)_p.
\end{equation*}
\end{lem}

From above and $\Delta_h^k (P, x, J)=0$ for all $P\in\Pi_k$ it readily follows that
\begin{equation}\label{E-equiv-omega}
E_k(f, J)_p \sim \omega_k(f, J)_p, \quad \forall f \in L^p(J), \;\; 1\le p\le \infty.
\end{equation}

The following useful property of modili of smoothness is well known:
\begin{equation}\label{omega-aver}
\omega_k(f, J)_p^p \sim \frac{1}{|J|}\int_0^{|J|}\int_J|\Delta_h^k (f, x, J)|^p dx dh,
\quad 1\le p<\infty.
\end{equation}
For proofs of the above claims and further details, see e.g. \cite[\S~7.1]{PP}.

We find useful the concept of {\it near best approximation}.
A polynomial $P_J(f) \in \Pi_k$ is said to be a polynomial of near best
$L^p(J)$-approximation to $f$ from $\Pi_k$ with constant $A\ge 1$ if
\begin{equation*}
\|f-P_J(f)\|_{L^p(J)} \le A E_k(f,J)_p.
\end{equation*}
Note that since $p \ge 1$ a near best
$L^p(J)$-approximation $P_J(f)$ (with an appropriate $A$) can be easily
realized by a linear projector, see below.
\begin{lem}{\cite[Chapter 12, Lemma 6.2]{DL}}\label{lem:near-best}
Suppose $1\le q< p\le\infty$ and $P_J$ is a polynomial of near best
$L^q(J)$-approximation to $f \in L^p(J)$ from $\Pi_k$.
Then $P_J$ is a polynomial of near best $L^p(J)$-approximation to $f$.
\end{lem}

\smallskip
\noindent
{\bf Local projectors onto polynomials.}
Given $1\le p\le \infty$ and a compact interval $J$,
we let $P_{J, p}: L^p(J)\to \Pi_k$ be a \emph{linear projector} such that
\begin{equation}\label{project}
\|f-P_{J,p}(f)\|_{L^p(J)} \le A E_k(f, J)_p,
\quad \forall f\in L^p(J),
\end{equation}
where the constant $A\ge 1$ is independent of $J$ and $f$.
Note that \eqref{project} implies $\|P_{J, p}\|_{L^p(J)\to L^p(J)}\le 1+A$, i.e. $P_{J, p}$ is a bounded linear operator.
$P_{J, p}$ can be defined via the averaged Taylor polynomial, see e.g. \cite[Section 4.1]{BS},
which is a linear operator for fixed $J, k$ that does not depend on $p$ (cf. Lemma~\ref{lem:near-best}).

\begin{rem}\label{rem:local-pol}
As is well know most of the above approximation results are valid in the wider range $p>0$.
The restriction $p\ge 1$ reflects our believe that
polynomial approximation of functions in $L^p$, $p<1$, is not natural.
There is no linear operator realization of polynomial near best approximation in $L^p$, $p<1$,
because there are no continuous linear functional in $L^p$, $p<1$.
Hence, the polynomial approximation in $L^p$, $p<1$, is \emph{nonconstructive}.
It should be replaced by approximation in the Hardy space $H^p$, $p<1$.
However, as will be seen later on the use of polynomial approximation in $L^p$ or $H^p$, $p<1$,
can be avoided completely in nonlinear spline approximation in $\BMO$.
\end{rem}

\subsection{Quasi-interpolant}\label{subsec:quasi}

We define the linear operator $T_m: \tilde\cS_m^k \to \cS_m^k$ by
\begin{equation*}
T_m(S):= \sum_{\Q\in\QQ_m} \aa_\Q(S)\varphi_\Q, \quad \forall S\in \tilde\cS_m^k,
\end{equation*}
where $\aa_\Q(S)$ are defined in \eqref{a-Q} with the additional requirement $\xi_\Q$ does not coincide with a knot of $\QQ_m$.
As observed in \S~\ref{subsec:splines} the identity
$T_m(S) = S$ for every $S\in \cS_m^k$ holds independently of the choice of $\xi_\Q$.
But for $S\in \tilde\cS_m^k \setminus \cS_m^k$ the value of $\aa_\Q(S)$ may depend on
the interval $(x_{m, \nu},x_{m, \nu+1})$, $\nu=j, \dots, j+k-1$, containing $\xi_\Q$.
In order to have well defined operator $T_m$ from now on for every $\Q\in\QQ_m$
we fixed the interval $I\in\II_m$, $I\subset \Q$, containing $\xi_\Q$ in its interior.

From \eqref{a-Q} it follows that
\begin{equation*}
|\aa_\Q(S)| \le c\|S\|_{L^\infty(\Q)}, \quad \forall S\in \tilde\cS_m^k,\quad \Q\in\QQ_m,
\end{equation*}
which easily leads to
\begin{equation}\label{Q-m-1}
\|T_m(S)\|_{L^p(I)} \le c\|S\|_{L^p(\Omega_I)},
\quad \forall S\in \tilde\cS_m^k, \; I\in\II_m, \; 1\le p\le \infty,
\end{equation}
where $\Omega_I$ is given in \eqref{def-Omega-I}.

We next extend $T_m$ to $L^p_\loc(\RR)$.
Let $P_{J, p}: L^p(J)\to \Pi_k$ be the projector from \eqref{project}.
Define
\begin{equation}\label{def-proj}
\sP_{m, p}(f):=\sum_{I\in\II_m} \ONE_I\cdot P_{I, p}(f),
\quad f\in L^p_\loc(\RR),
\end{equation}
(the precise values of $\sP_{m, p}(f)$ at the knots of $\II_m$ are not of importance for this study) and set
\begin{equation}\label{def-proj-2}
T_{m, p}(f):= T_m(\sP_{m, p}(f)), \quad f\in L^p_\loc(\RR).
\end{equation}

\begin{lem}\label{lem:quasi-int}
If $f\in L^p_\loc(\RR)$, $1\le p\le \infty$, $m\in\ZZ$, then
\begin{equation*}
\|f- T_{m, p}(f)\|_{L^p(I)} \le cE_k(f, \Omega_I)_p, \quad \forall I\in \II_m.
\end{equation*}
\end{lem}
\proof
Let $R\in\Pi_k$ be such that $\|f-R\|_{L^p(\Omega_I)} \le cE_k(f, \Omega_I)_p$.
Then
\begin{align*}
\|f- T_{m, p}(f)\|_{L^p(I)} &= \|f- T_m(\sP_{m, p}(f))\|_{L^p(I)}
\\
& \le \|f-R\|_{L^p(I)} + \|R-T_m(\sP_{m, p}(f))\|_{L^p(I)}
\\
& = \|f-R\|_{L^p(I)} + \|T_m[\sP_{m, p}(R-f)]\|_{L^p(I)}
\\
& \le \|f-R\|_{L^p(\Omega_I)} + c\|\sP_{m, p}(R-f)\|_{L^p(\Omega_I)}
\\
& \le c\|f-R\|_{L^p(\Omega_I)}\le cE_k(f, \Omega_I)_p.
\end{align*}
Here we used \eqref{Q-m-1} and the boundedness of $\sP_{m, p}$.
$\hfill \square$

\subsection{Embedding of sequence type B-spline spaces into \boldmath{$\BMO$}}\label{subsec:embed-in-BMO}

The sequence $\ell^\tau = \ell^\tau (\QQ)$ space indexed by $\QQ$, $0<\tau\le\infty$,
is the set of all sequence $\{\cc_\Q\}_{\Q\in\QQ}$ of complex numbers such that
\begin{equation*}
\|\{\cc_\Q\}\|_{\ell^\tau}
:= \Big(\sum_{\Q\in\QQ} |\cc_\Q|^\tau\Big)^{1/\tau}
<\infty.
\end{equation*}

The embedding of the Besov spaces of interest to us in $\BMO$
will play a crucial role in this article.
This embedding is in essence contained in the following

\begin{thm}\label{thm:embed-B-BMO}
Let $\{\cc_\Q\}_{\Q\in\QQ}\in\ell^\tau (\QQ)$,  $0<\tau<\infty$.

$(a)$ If $\tau>1$,
then $\sum_{\Q\in\QQ} \cc_\Q\varphi_\Q$ converges unconditionally in $\BMO(\RR)$
and
\begin{equation}\label{embed-B-BMO}
\Big\|\sum_{\Q\in\QQ} \cc_\Q\varphi_\Q\Big\|_{\BMO} \le c\|\{\cc_\Q\}\|_{\ell^\tau}.
\end{equation}
Consequently, $\sum_{\Q\in\QQ} \cc_\Q\varphi_\Q \in \VMO(\RR)$.

$(b)$ If $0<\tau\le 1$, then obviously
$\sum_{\Q\in\QQ} \cc_\Q\varphi_\Q$ converges absolutely and unconditionally in $L^\infty(\RR)$
to a~continuous function and
\begin{equation*}
\Big\|\sum_{\Q\in\QQ} \cc_\Q\varphi_\Q\Big\|_{\BMO}
\le \Big\|\sum_{\Q\in\QQ} |\cc_\Q\varphi_\Q|\Big\|_{\infty} \le c\|\{\cc_\Q\}\|_{\ell^\tau}.
\end{equation*}
\end{thm}

The proof of this theorem depends on the following

\begin{lem}\label{lem:embed-BMO}
Let $0<p, \tau<\infty$.
Then for any sequence $\{\cc_\Q\}_{\Q\in\QQ}$ and any compact interval $J\subset \RR$
we have
\begin{equation}\label{embed-BMO}
\Big\|\frac{1}{|J|^{1/p}}\sum_{\Q\in\QQ, \Q\subset J}|\cc_\Q\varphi_\Q|\Big\|_p
\le c\Big(\sum_{\Q\in\QQ, \Q\subset J}|\cc_\Q|^\tau\Big)^{1/\tau}.
\end{equation}
\end{lem}

This lemma is a consequence of the following well known embedding result (see e.g. \cite[Theorem~3.3]{KP}).

\begin{prop}\label{prop:embed}
If $0<\tau < p<\infty$, then for any sequence $\{\cc_\Q\}_{\Q\in\QQ}$ of complex number
one has
\begin{equation*}
\Big\|\sum_{\Q\in\QQ} |\cc_\Q\varphi_\Q| \Big\|_p
\le c\Big(\sum_{\Q\in\QQ} \|\cc_\Q\varphi_\Q \|_p^\tau\Big)^{1/\tau},
\end{equation*}
where the constant $c>0$ depends only on $p$, $\tau$,
and the parameter $\lambda$ from condition $(c)$ on $\II$.
\end{prop}

To streamline our presentation we defer the proofs of Lemma~\ref{lem:embed-BMO} and Theorem~\ref{thm:embed-B-BMO}
to \S\ref{subsec:appendix1} in the appendix.

\section{Homogeneous Besov spaces}\label{sec:Besov-spaces}

In this section, we introduce and discuss the Besov spaces $\BB^{\alpha, k}_\tau$
that will be used for characterization of nonlinear $n$-term spline approximation in $\BMO$.
For the theory of Besov spaces we refer the reader to \cite{Peetre, T1, FJ1, FJ2, FJW}.

Throughout the section we assume that
\begin{equation*}
\alpha>0, \; k\ge 2, \quad\hbox{and}\quad \tau:=1/\alpha.
\end{equation*}

\subsection{The homogeneous Besov spece \boldmath{$\BB^{\alpha, k}_\tau$} in the case \boldmath{$\tau\ge 1$}}

\begin{defn}\label{def:B-1}
The homogeneous Besov space $\BB^{\alpha, k}_\tau=\BB^{\alpha, k}_\tau(\RR)$, $1\le \tau <\infty$ $(0<\alpha\le 1)$, is defined as the collection of all
functions $f\in L^\tau_\loc(\RR)$ such that
$\Delta^k_hf \in L^\tau(\RR)$ for all $h\in\RR$ and
\begin{equation}\label{def-B-norm-modulus}
\|f\|_{\BB^{\alpha, k}_\tau}
:= \Big(\int_0^\infty [t^{-\alpha}\omega_k(f, t)_\tau]^\tau\frac{dt}{t}\Big)^{1/\tau} <\infty,
\end{equation}
where
\begin{equation*}
\omega_k(f, t)_\tau := \sup_{|h|\le t} \|\Delta_h^k f\|_{\tau},\quad
\Delta_h^k f( x)=\sum_{j=0}^k (-1)^{k+j} \binom{k}{j} f(x+jh).
\end{equation*}
\end{defn}
Notice the different definition and notation of finite differences \eqref{delta-f} and  moduli \eqref{omega-J} on compact interval.
Observe that
$\|f+P\|_{\BB^{\alpha, k}_\tau} = \|f\|_{\BB^{\alpha, k}_\tau}$
for each polynomial $P\in\Pi_k$.

From the properties of $\omega_k(f, t)_\tau$ it readily follows that
\begin{equation}\label{discr-B-mogulus}
\|f\|_{\BB^{\alpha, k}_\tau}
\sim \Big(\sum_{\nu\in\ZZ}\big(2^{\alpha \nu}\omega_k(f, 2^{-\nu})_\tau\big)^\tau\Big)^{1/\tau}.
\end{equation}
It is also easy to see that
\begin{equation}\label{B-tau}
\|f\|_{\BB^{\alpha, k}_\tau}\sim  \Big(\sum_{I \in \II} (|I|^{-\alpha} \omega_k(f, \Om_I)_\tau)^\tau\Big)^{1/\tau}
\sim \Big(\sum_{I \in \II} (|I|^{-\alpha} E_k(f, \Om_I)_\tau)^\tau\Big)^{1/\tau}.
\end{equation}
Recall that $E_k(f, \Om_I)_\tau\sim\omega_k(f, \Om_I)_\tau$, see \eqref{E-equiv-omega}.
For reader's convenience we give a simple proof of equivalence \eqref{B-tau} in \S\ref{subsec:appendix2} in the appendix.

Subsequently, in \thmref{thm:B-equiv-norms} we shall show the following key equivalence: For any $1\le q<\infty$, $\tau\ge 1$,
\begin{equation}\label{B-q}
\|f\|_{\BB^{\alpha, k}_\tau}
\sim \Big(\sum_{I \in \II} \big(|I|^{-\frac{1}{q}} \omega_k(f, \Om_I)_q\big)^\tau\Big)^{1/\tau}
\sim \Big(\sum_{I \in \II} \big(|I|^{-\frac{1}{q}} E_k(f, \Om_I)_q\big)^\tau\Big)^{1/\tau}.
\end{equation}

\subsection{The Besov space \boldmath{$\BB^{\alpha, k}_\tau$} in the general case when \boldmath{$0<\tau <\infty$}}

As was alluded to in Definition~\ref{rem:local-pol} to us polynomial approximation in $L^\tau$, $\tau<1$,
is not normal and should be avoided.
Furthermore, we think that even when $\tau\ge 1$ the Besov spaces $\BB^{\alpha, k}_\tau$ are most naturally defined
via local polynomial approximation in $L^q$ with $q\ge 1$.
The equivalence \eqref{B-q} is a motivation for making the following

\begin{defn}\label{def:B-2}
Let $\II$ be a regular multilevel partition of $\RR$ $($see \S\ref{subsec:nested_partitions}$)$.
The Besov space $\BB^{\alpha, k}_{\tau}(E, q)$, $1\le q<\infty$, $0<\tau<\infty$, is defined as the collection of all
functions $f\in L^q_\loc(\RR)$ such that
\begin{equation}\label{def-B-q}
\|f\|_{\BB^{\alpha, k}_\tau(E, q)}
:= \Big(\sum_{I \in \II} \big(|I|^{-\frac{1}{q}} \omega_k(f, \Om_I)_q\big)^\tau\Big)^{1/\tau}
\sim \Big(\sum_{I \in \II} \big(|I|^{-\frac{1}{q}} E_k(f, \Om_I)_q\big)^\tau\Big)^{1/\tau}
\end{equation}
is finite.
Here
$\omega_k(f, \Om_I)_q$ and $E_k(f, \Om_I)_q$ are defined in \eqref{omega-J} and \eqref{Ek-J}
with $\Omega_I$ from \eqref{def-Omega-I}.
\end{defn}

\begin{lem}\label{lem:independ}
The Besov space $\BB^{\alpha, k}_{\tau}(E, q)$ introduced above
is independent of the particular selection of the multilevel partition $\II$ of $\RR$ being used.
\end{lem}
\proof
Let $\II'$ be another  multilevel partition of $\RR$ with the properties of $\II$.
It is readily seen that for every level $\II_m'$ of $\II'$ there exists
a level $\II_n$ of $\II$ with these properties:
(a) The intervals in $\II_m'$ and $\II_n$ are comparable in length, and
(b) For each interval $I'\in\II'_m$ there exists an interval $I\in\II_n$
such that $\Omega_{I'}\subset\Omega_I$.
For example, condition (b) is satisfied whenever $\max_{I'\in\II_m'}|I'|\le 3/2 \min_{I\in\II_n}|I|$.
Hence,
\begin{equation*}
\sum_{I' \in \II'_m} \big(|I'|^{-\frac{1}{q}} \omega_k(f, \Om_{I'})_q\big)^\tau
\le c\sum_{I \in \II_n} \big(|I|^{-\frac{1}{q}} \omega_k(f, \Om_I)_q\big)^\tau.
\end{equation*}
Clearly, each level $\II_n$ of $\II$ can serve in this capacity for only uniformly bounded number
of levels from $\II'$.
The claim follows.
$\hfill \square$

In light of \eqref{B-tau} the spaces $\BB^{\alpha, k}_{\tau}$ and $\BB^{\alpha, k}_{\tau}(E, \tau)$
are the same with equivalent norms when $\tau\ge 1$.
Furthermore, as will be shown in \thmref{thm:B-equiv-norms} \eqref{B-q} is valid and hence
the spaces $\BB^{\alpha, k}_{\tau}(E, q)$ are the same space for all $1\le q<\infty$ with equivalent norms when $\tau\ge 1$.
The same theorem extends the equivalence of these spaces for $0<\tau<\infty$.
To achieve this we need some preparation.

\subsection{Norms via projectors}

We define
\begin{equation}\label{def:qt}
\qq_{m, q} := T_{m,q} - T_{m-1,q}
\;\;\mbox{for}\;\; m \in \ZZ,
\end{equation}
where $T_{m,q}$ is defined in \eqref{def-proj-2}.
For a given function
$f \in L^q_\loc(\RR)$, $1\le q < \infty$, clearly
$\qq_{m,q}(f) \in \cS_{m}^k$
and we define uniquely the sequence
$\{b_{\Q, q}(f)\}_{\Q \in \QQ_m}$ by
\begin{equation}\label{def:b-tm}
\qq_{m,q}(f)
=: \sum_{\Q \in \QQ_m} b_{\Q, q}(f) \varphi_\Q.
\end{equation}
Since the near-best approximant in \eqref{project} is realized as a linear operator, then, evidently,
$\{b_{\Q, q}(\cdot)\}$
are linear functionals.

We introduce the following (quasi-)norms for functions $f\in L^q_{\loc}(\RR)$, $1~\le~q<\infty$, $0<\tau<\infty$:
\begin{equation*}
\|f\|_{\BB^{\alpha, k}_\tau(Q, q)}
:= \Big( \sum_{\Q \in \QQ} (|\Q|^{-1/q} \|b_{\Q,q}(f) \varphi_\Q \|_q)^\tau \Big)^{1/\tau}.
\end{equation*}
By \eqref{norm-spline} we have
\begin{equation}\label{NQ=b-normt}
\|f\|_{\BB^{\alpha, k}_\tau(Q, q)}
\sim \Big(\sum_{I \in \II}(|I|^{-1/q}\|\qq_{m, q}(f)\|_{L^q(I)})^\tau \Big)^{1/\tau},
\end{equation}
and, since $\|\varphi_\Q\|_p \sim |\Q|^{1/p}$,
\begin{equation}\label{NQ=b-norm1}
\|f\|_{\BB^{\alpha, k}_\tau(Q,q)} \sim \Big(\sum_{\Q \in \QQ} | b_{\Q,q}(f)|^\tau \Big)^{1/\tau}.
\end{equation}

\begin{lem}\label{lem:BQ-BE}
If $f\in \BB^{\alpha, k}_\tau(E,q)$, $1\le q<\infty$, then
\begin{equation}\label{BQ-BE}
\|f\|_{\BB^{\alpha, k}_\tau(Q,q)} \le c \|f\|_{\BB^{\alpha, k}_\tau(E,q)}.
\end{equation}
\end{lem}
\proof
Let $f\in \BB^{\alpha, k}_\tau(E,q)$.
If $I\in\II_j$ and $J\in \II_{j-1}$ is the unique parent of $I$ ($I\subset J$),
then by Lemma~\ref{lem:quasi-int}
\begin{align}\label{est-qj}
\|\qq_{j, q}(f)\|_{L^q(I)} &\le c\|f-T_{j, q}(f)\|_{L^q(I)} + c\|f-T_{j-1, q}(f)\|_{L^q(J)}
\\
&\le c E_k(f, \Omega_I)_q + c E_k(f, \Omega_J)_q,
\quad 1\le q<\infty. \nonumber
\end{align}
This, \eqref{NQ=b-normt} and \eqref{def-B-q} imply \eqref{BQ-BE}.
$\hfill \square$

As will be shown later $\|\cdot\|_{\BB^{\alpha, k}_\tau(Q,q)}$ is another
equivalent norm in $\BB^{\alpha, k}_\tau(E,q)$.

\subsection{Decomposition of \boldmath{$\BB^{\alpha, k}_\tau(E, q)$} and embedding in $\VMO$}

Our next step is to derive a representation of the functions in $\BB^{\alpha, k}_\tau(E, q)$
via the quasi-intrepolant from \eqref{def-proj-2}.
We first show that the Besov space $\BB^{\alpha, k}_\tau(E, q)$
is embedded in $\BMO$ modulo polynomials of degree $k-1$.

We define the $\BMO$ type space $\BMO^{q,k}(\RR)$, $1\le q<\infty$, as the set of all functions $f\in L^q_{\loc}(\RR)$ such that
\begin{equation*}
\|f\|_{\BMO^{q,k}}:=\sup_{I\in\II}|I|^{-1/q}\omega_k(f,\Omega_I)_q\sim\sup_{I\in\II}|I|^{-1/q}E_k(f,\Omega_I)_q<\infty.
\end{equation*}

\begin{prop}\label{prop:BMO-embed-3}
For any $f\in \BMO^{q,k}(\RR)$, $1\le q<\infty$,
there exists a polynomial $P\in\Pi_k$ such that $f-P\in\BMO$ and
\begin{equation}\label{B-embed-BMO-3}
\|f-P\|_{\BMO} \le c \|f\|_{\BMO^{q,k}}.
\end{equation}
\end{prop}
\proof
If $\|f\|_{\BMO^{q,k}}=0$ then $f$ coincides with its polynomial of best approximation on every of the over-lapping intervals $\Omega_I$
and \eqref{B-embed-BMO-3} follows trivially.

Let $\|f\|_{\BMO^{q,k}}>0$.
Denote
$J_\nu:= [-2^\nu, 2^\nu]$, $\nu\in \NN_0$.
Evidently, for each $\nu\in\NN_0$ there exists an interval $I_\nu\in \II$
such that $J_\nu\subset \Omega_{I_\nu}$ and $|J_\nu|\sim |I_\nu|$.
In light of Whitney's theorem (Lemma~\ref{lem:whitney}), there exists a polynomial $P_\nu\in\Pi_k$ such that
\begin{align}\label{embed-11}
\frac{1}{|J_\nu|}\int_{J_\nu}|f(x)-P_\nu(x)|^q dx
\le c|J_\nu|^{-1}\omega_k(f, J_\nu)_q^q
\le c\|f\|_{\BMO^{q,k}}^q.
\end{align}
Denote
$\Upsilon_\nu:=P_\nu-P_{\nu-1}$.
Since $\Upsilon_\nu\in\Pi_k$ we have for $x\in J_n$, $n\ge 0$, and $\nu\ge n$
\begin{equation*}
|\Upsilon_\nu(x)-\Upsilon_\nu(0)|
\le |J_n|\|\Upsilon_\nu'\|_{L^\infty(J_n)}
\le |J_n|\|\Upsilon_\nu'\|_{L^\infty(J_\nu)}
\le c|J_n||J_\nu|^{-1}\|\Upsilon_\nu\|_{L^\infty(J_\nu)}
\end{equation*}
and using \lemref{lem:poly-norms} and \eqref{embed-11}
\begin{align*}
\|\Upsilon_\nu\|_{L^\infty(J_\nu)}
&\le |J_\nu|^{-1/q}\|\Upsilon_\nu\|_{L^q(J_\nu)}\le |J_\nu|^{-1/q}\|\Upsilon_\nu\|_{L^q(J_{\nu-1})}
\\
&\le c|I_\nu|^{-1/q}\omega_k(f, \Omega_{I_\nu})_q
+c|I_{\nu-1}|^{-1/q}\omega_k(f, \Omega_{I_{\nu-1}})_q
\le c\|f\|_{\BMO^{q,k}}.
\end{align*}
Then for any $n, m\in\NN$, $n<m$, we have
\begin{align}\label{L-L}
\sum_{\nu=n+1}^m \|\Upsilon_\nu-\Upsilon_\nu(0)\|_{L^\infty(J_n)}
&\le c|J_n|\sum_{\nu=n+1}^m |J_\nu|^{-1}\|\Upsilon_\nu\|_{L^\infty(J_\nu)}
\\
&\le c\|f\|_{\BMO^{q,k}}. \nonumber
\end{align}
From above with $n=0$ it follows that the series
$\sum_{\nu=1}^\infty (\Upsilon_\nu-\Upsilon_\nu(0))$
converges uniformly on $J_0=[-1, 1]$ to some polynomial in $\Pi_k$.
Hence, there exists a polynomial $P\in\Pi_k$ such that
\begin{equation*} 
\|P_m-P_m(0)-(P-P(0))\|_{L^\infty(J_0)} \to 0
\quad\hbox{as}\quad m\to\infty.
\end{equation*}
From this and \eqref{L-L} it follows that for any (fixed) $n\in\NN$
\begin{equation}\label{P-m-Pn}
\|P_m-P_m(0)-(P-P(0))\|_{L^\infty(J_n)} \to 0
\quad\hbox{as}\quad m\to\infty.
\end{equation}

We shall next show that \eqref{B-embed-BMO-3} holds with the polynomial $P$ from above.
Let $J$ be an arbitrary compact interval.
Then there exists an interval $\tilde{I}_0\in\II$ such that $J\subset \Omega_{\tilde{I}_0}$ and $|\tilde{I}_0|\sim |J|$.
Let $n\in\NN$ be the minimal positive integer such that $J\subset J_n$.
Let $\{\tilde{I}_j\}_{j=0}^\ell$ be intervals from consecutive levels of $\II$ such that
$\tilde{I}_0 \subset \tilde{I}_1\subset \cdots\subset \tilde{I}_\ell$, $\tilde{I}_j$ is a parent of $\tilde{I}_{j-1}$,
$\tilde{I}_\ell \cap J_n\ne \emptyset$, and $|\tilde{I}_\ell|\sim |J_n|$.

As $n\in \NN$ is already fixed we choose $m>n$ so that
\begin{equation}\label{Pm-P-B}
\|P_m-P_m(0)-(P-P(0))\|_{L^\infty(J_n)}< \|f\|_{\BMO^{q,k}}.
\end{equation}
This is possible because of \eqref{P-m-Pn}.

Just as in \eqref{embed-11}, applying Whitney's theorem (Lemma~\ref{lem:whitney})
there exist polynomials $\tilde{P}_{j}\in\Pi_k$, $j=0,1,\dots,\ell$, such that
\begin{equation}\label{whitney-2}
\frac{1}{|\tilde{I}_j|}\int_{\Omega_{\tilde{I}_j}}|f(x)-\tilde{P}_{j}(x)|^q dx
\le c|\tilde{I}_j|^{-1}\omega_k(f, \Omega_{\tilde{I}_j})_q^q
\le c\|f\|_{\BMO^{q,k}}^q.
\end{equation}

At this point we select several constants.
We choose $\tilde{c}:=\tilde{P}_0(y)-\tilde{P}_\ell(y)$, where $y\in J$ is fixed,
$c^\star := P_n(0)-P_m(0)$, and
$c^{\star\star}:= P_m(0)-P(0)$.
We also set $c^\diamond=\tilde{c}+c^\star+c^{\star\star}$.

Using the above polynomials and constants we get
\begin{align*}
\frac{1}{|J|}\int_J &|f(x)-P(x)-c^\diamond|^qdx
\le \frac{c}{|J|}\int_J|f(x)-\tilde{P}_0(x)|^qdx
\\
&+ c\|\tilde{P}_0-\tilde{P}_\ell-\tilde{c}\|_{L^\infty(J)}^q
+ c\|\tilde{P}_\ell-P_n\|_{L^\infty(J)}^q
+ c\|P_n-P_m-c^\star\|_{L^\infty(J)}^q
\\
&+ c\|P_m-P-c^{\star\star}\|_{L^\infty(J)}^q
=: W_1+W_2+W_3+W_4+W_5.
\end{align*}
To estimate $W_1$ we use \eqref{whitney-2} and obtain
\begin{equation*}
W_1 \le \frac{1}{|\tilde{I}_0|}\int_{\Omega_{\tilde{I}_0}}|f(x)-\tilde{P}_{j}(x)|^q dx
\le  c|\tilde{I}_0|^{-1}\omega_k(f, \Omega_{\tilde{I}_0})_q^q
\le c\|f\|_{\BMO^{q,k}}^q.
\end{equation*}
We proceed just as in \eqref{L-L} to obtain
\begin{align*}
W_2=c\|\tilde{P}_0-\tilde{P}_\ell-\tilde{c}\|_{L^\infty(J)}^q
\le c\|f\|_{\BMO^{q,k}}^q.
\end{align*}
We now estimate $W_3$.
Observe that because $\tilde{I}_\ell \cap J_n\ne \emptyset$ and $|\tilde{I}_\ell|\sim |J_n|$
we have $|\Omega_{\tilde{I}_\ell}\cap J_n|\sim |\Omega_{\tilde{I}_\ell}|\sim |J_n|$.
We use this, \eqref{embed-11}, and \eqref{whitney-2} to obtain
\begin{align*}
W_3 &\le c\|\tilde{P}_\ell-P_n\|_{L^\infty(J_n\cap\Omega_{\tilde{I}_\ell})}^q
\le c|J_n\cap\Omega_{\tilde{I}_\ell}|^{-1}\|\tilde{P}_\ell-P_n\|_{L^q(J_n\cap\Omega_{\tilde{I}_\ell})}^q
\\
& \le c|J_n|^{-1}\|f-P_n\|_{L^q(J_n)}^q + c|\Omega_{\tilde{I}_\ell}|^{-1}\|f-\tilde{P}_\ell\|_{L^q(\Omega_{\tilde{I}_\ell})}^q
\le c\|f\|_{\BMO^{q,k}}^q.
\end{align*}
To estimate $W_4$ we use \eqref{L-L} and obtain
$
W_4 \le c\|f\|_{\BMO^{q,k}}^q.
$
We also have from \eqref{Pm-P-B} that
$W_5\le c\|f\|_{\BMO^{q,k}}^q$.

Putting the above estimates together we obtain that for any interval $J$
there is a constant $c^\diamond$ such that
\begin{equation*}
\frac{1}{|J|}\int_J |f(x)-P(x)-c^\diamond|^qdx
\le c\|f\|_{\BMO^{q,k}}^q
\end{equation*}
and \eqref{B-embed-BMO-3} follows in view of \eqref{equiv-norms-BMO}.
$\hfill \square$

From the trivial embedding $\BB^{\alpha, k}_\tau(E, q)\subset\BMO^{q,k}$,
which is a consequence of $\|\cdot\|_{\ell^\infty}\le\|\cdot\|_{\ell^\tau}$,
and \propref{prop:BMO-embed-3} we obtain
\begin{prop}\label{prop:BMO-embed}
For any $f\in \BB^{\alpha, k}_\tau(E, q)$, $1\le q<\infty$, $0<\tau<\infty$,
there exists a polynomial $P\in\Pi_k$ such that $f-P\in\BMO$ and
\begin{equation*}
\|f-P\|_{\BMO} \le c \|f\|_{\BB^{\alpha, k}_\tau(E, q)}.
\end{equation*}
\end{prop}

The following decomposition result will play a central role in our theory of Besov spaces.

\begin{thm}\label{thm:rep-B}
For any $f\in \BB^{\alpha, k}_\tau(E, q)$, $1\le q<\infty$, $0<\tau<\infty$,
there exists a polynomial $P\in\Pi_k$ such that
$f-P\in \VMO$ and
\begin{equation}\label{rep-B}
f-P =\sum_{m\in\ZZ}\sum_{\Q\in\QQ_m} b_{Q,q}(f)\varphi_\Q,
\end{equation}
where the convergence is unconditional in $\BMO$.
Here the coefficients $\{b_{Q,q}(f)\}$ are from \eqref{def:b-tm}.
Furthermore,
\begin{equation}\label{BE-embed-BMO}
\Big(\sum_{\Q \in \QQ} | b_{\Q,q}(f)|^\tau \Big)^{1/\tau}
\le c \|f\|_{\BB^{\alpha, k}_\tau(E, q)}.
\end{equation}
\end{thm}

We divert the long and tedious proof of this theorem to \S\ref{subsec:appendix3} in the appendix.

\subsection{Norm in \boldmath{$\BB^{\alpha, k}_\tau$} via B-splines}

Theorems \ref{thm:embed-B-BMO} and \ref{thm:rep-B} are the motivation for the following

\begin{defn}\label{def:Besov-B-spline}
The Besov space $\BB^{\alpha, k}_\tau(\FF)$ is defined as the collection of all
functions $f$ on $\RR$ for which there exists a polynomial $P\in\Pi_k$  such that
$f-P\in \VMO$ and $f-P$
can be represented in the form
$f-P=\sum_{\Q\in\QQ} \cc_\Q\varphi_\Q$,
where the convergence is unconditional in $\BMO$,
and $\sum_{\Q \in \QQ}| \cc_\Q |^\tau<\infty$.
The norm in $\BB^{\alpha,k}_\tau(\FF)$ is defined by
\begin{equation*}
\|f\|_{\BB^{\alpha,k}_\tau(\FF)}
:=\inf\Big\{\Big( \sum_{\Q \in \QQ}| \cc_\Q|^\tau \Big)^{1/\tau} : f-P= \sum_{\Q \in \QQ} \cc_\Q \varphi_\Q\Big\}.
\end{equation*}
\end{defn}

Observe that due to $\|\varphi_\Q\|_q \sim |\Q|^{1/q}$, $0<q\le \infty$, we have

\begin{equation*}
\|f\|_{\BB^{\alpha, k}_\tau(\FF)}
\sim
\inf\Big\{\Big( \sum_{\Q \in \QQ}\big(|\Q|^{-1/q}\| \cc_\Q \varphi_\Q \|_q\big)^\tau \Big)^{1/\tau} : f-P= \sum_{\Q \in \QQ} \cc_\Q \varphi_\Q\Big\}.
\end{equation*}

\subsection{Equivalent Besov norms}

\begin{thm}\label{thm:B-equiv-norms}
The homogeneous Besov spaces $\BB^{\alpha, k}_\tau(\FF)$, $\BB^{\alpha, k}_\tau(E, q)$ and $\BB^{\alpha, k}_\tau(Q, q)$ for all $q\in [1, \infty)$
are the same with equivalent norms:
\begin{equation*}
\|f\|_{\BB^{\alpha, k}_\tau(\FF)} \sim \|f\|_{\BB^{\alpha, k}_\tau(E, q)} \sim \|f\|_{\BB^{\alpha, k}_\tau(Q, q)}
\end{equation*}
with constants of equivalence depending only on $\alpha$, $k$, $q$,
and the parameters of $\II$.
\end{thm}
\proof
(a) We first show that if $f\in \BB^{\alpha, k}_\tau(\FF)$, then $f\in \BB^{\alpha, k}_\tau(E, q)$
for $q\in [1, \infty)$ 
and
\begin{equation}\label{B-E-phi}
\|f\|_{\BB^{\alpha, k}_\tau(E, q)} \le c \|f\|_{\BB^{\alpha, k}_\tau(\FF)}.
\end{equation}

Note that H\"{o}lder's inequality implies for any compact interval $J$
\begin{equation*}
|J|^{-1/p}E_k(f, J)_p \le |J|^{-1/q}E_k(f, J)_q,\quad 1\le p<q\le\infty,\quad f\in L^q(J),
\end{equation*}
and hence
\begin{equation*}
\|f\|_{\BB^{\alpha, k}_\tau(E, p)} \le c \|f\|_{\BB^{\alpha, k}_\tau(E, q)},
\quad\hbox{if}\quad 1\le p \le q <\infty.
\end{equation*}
Therefore, it is sufficient to prove \eqref{B-E-phi}
in the (most unfavorable) case when $q>\max\{1, \tau\}$.

Assume $f\in \BB^{\alpha, k}_\tau(\FF)$. Then by Definition~\ref{def:Besov-B-spline}
$f$ can be represented in the form
\begin{equation*}
f-P= \sum_{\Q\in\QQ}\cc_\Q\varphi_\Q
\end{equation*}
with unconditional convergence in $\BMO$, where
$P\in\Pi_k$, $f-P\in\VMO$ and the coefficients $\{\cc_\Q\}$
satisfy
$\big(\sum_{\Q\in\QQ}|\cc_\Q|^\tau\big)^{1/\tau}
\le 2\|f\|_{\BB^{\alpha, k}_\tau(\FF)}$.

Denote by $\ell(I)$ the level of $I$ ($\ell(I)=m$ if $I\in\II_m$)
and, similarly, by $\ell(\Q)$ the level of $\Q$.
Fix $I\in\II$.
To estimate $\omega_k(f, \Omega_I)_q$ we split the representation of $f-P$ into two:
\begin{equation*}
f-P= \sum_{j\ge \ell(I)}\sum_{\Q\in\QQ_j}\cc_\Q\varphi_\Q + \sum_{j<\ell(I)}\sum_{\Q\in\QQ_j}\cc_\Q\varphi_\Q
=: G_I+H_I.
\end{equation*}
We use Proposition~\ref{prop:embed} to obtain
\begin{align}\label{est-omega-g}
\omega_k(G_I, \Omega_I)_q \le c\|G_I\|_{L^q(\Omega_I)}
&\le c\Big(\sum_{\substack{\ell(\Q)\ge \ell(I),\\ \Q\cap \Omega_I\ne \emptyset}} \|\cc_\Q\varphi_\Q\|_q^\tau\Big)^{1/\tau}
\\
&\le c\Big(\sum_{\substack{\ell(\Q)\ge \ell(I),\\ \Q\cap \Omega_I\ne \emptyset}}
|\Q|^{\tau/q}|\cc_\Q|^\tau\Big)^{1/\tau}. \nonumber
\end{align}
To estimate $\omega_k(H_I, \Omega_I)_q$ we use that
$\|\Delta^k_h\varphi_\Q\|_\infty \le c(|h|/|\Q|)^k\le c|h|/|\Q|$ (for $|h|\le |\Omega_I|/k$ and $\Q\in\QQ_j$ with $j<\ell(I)$)
and
$\Delta^k_h\varphi_\Q(x)=0$ if $[x, x+ kh]\cap V_\Q=\emptyset$,
where $V_Q$ is the set of all knots of $\varphi_\Q$.
We obtain
\begin{align}\label{est-omega-h}
\omega_k(H_I, \Omega_I)_q
&\le \sum_{j<\ell(I)}\sum_{\Q\in\QQ_j}\omega_k(\cc_\Q\varphi_\Q, \Omega_I)_q
\\
&\le c\sum_{j<\ell(I)}\sum_{\substack{\Q\in\QQ_j,\\ V_\Q\cap\Omega_I\ne\emptyset}}
\frac{|I|^{1+1/q}}{|\Q|}|\cc_\Q|. \nonumber
\end{align}
Using \eqref{def-B-q} we have
\begin{equation}\label{B-omega-omega}
\|f\|_{\BB^{\alpha, k}_\tau(E,q)}^\tau
\le c\sum_{I\in\II} |I|^{-\tau/q}\omega_k(G_I, \Omega_I)_q^\tau
+ c\sum_{I\in\II} |I|^{-\tau/q}\omega_k(H_I, \Omega_I)_q^\tau
=: S_1+S_2.
\end{equation}
To estimate $S_1$ we use \eqref{est-omega-g} and obtain
\begin{align*}
S_1 &\le \sum_{I\in\II} |I|^{-\tau/q}
\sum_{\substack{\ell(\Q)\ge \ell(I),\\ \Q\cap \Omega_I\ne \emptyset}}
|\Q|^{\tau/q}|\cc_\Q|^\tau
\\
&= c\sum_{\Q\in\QQ} |\cc_\Q|^\tau
\sum_{\substack{I\in\II,\\ \ell(I)\le \ell(\Q),\\ \Omega_I\cap \Q\ne \emptyset}} (|\Q|/|I|)^{\tau/q}
\le c\sum_{\Q\in\QQ} |\cc_\Q|^\tau
\sum_{\nu\ge 0}\rho^{\nu\tau/q}
\le c\|f\|_{\BB^{\alpha, k}_\tau(\FF)}^\tau.
\end{align*}
Here $\rho<1$ is from \eqref{r-rho} and we used the nested structure of the partition $\II$;
we switched once the order of summation.

We now estimate $S_2$. If $\tau\ge 1$ using \eqref{est-omega-h} we get
\begin{align*}
S_2 &\le c\sum_{I\in\II} |I|^{-\tau/q}
\Big(\sum_{j<\ell(I)}\sum_{\substack{\Q\in\QQ_j,\\ V_\Q\cap\Omega_I\ne\emptyset}}
\frac{|I|^{1+1/q}}{|\Q|}|\cc_\Q|\Big)^\tau
\\
&= c\sum_{I\in\II}
\Big(\sum_{j<\ell(I)}\sum_{\substack{\Q\in\QQ_j,\\ V_\Q\cap\Omega_I\ne\emptyset}}
(|I|/|\Q|)|\cc_\Q|\Big)^\tau
\\
& \le c\sum_{I\in\II}
\Big(\sum_{j<\ell(I)}\sum_{\substack{\Q\in\QQ_j,\\ V_\Q\cap\Omega_I\ne\emptyset}}
|I|/|\Q| \Big)^{\tau/\tau'}
\sum_{j<\ell(I)}\sum_{\substack{\Q\in\QQ_j,\\ V_\Q\cap\Omega_I\ne\emptyset}}
(|I|/|\Q|)|\cc_\Q|^\tau
\\
& \le c\sum_{I\in\II}
\sum_{j<\ell(I)}\sum_{\substack{\Q\in\QQ_j,\\ V_\Q\cap\Omega_I\ne\emptyset}}
(|I|/|\Q|)|\cc_\Q|^\tau.
\end{align*}
For the former inequality above we used H\"{o}lder's inequality.
For the last inequality we used \eqref{r-rho} as in the estimate of $S_1$.
Switching the order of summation in the last sums above we get
\begin{align*}
S_2\le c\sum_{\Q\in\QQ} |\cc_\Q|^\tau
\sum_{\substack{I\in\II,\\ \ell(I) > \ell(\Q),\\ \Omega_I\cap V_\Q\ne\emptyset}} |I|/|\Q|
\le c\sum_{\Q\in\QQ} |\cc_\Q|^\tau
\le c\|f\|_{\BB^{\alpha, k}_\tau(\FF)}^\tau.
\end{align*}
Here we used the simple fact that for any point $y\in\RR$ and $J\in \II$ we have
$$
\sum_{\substack{I\in\II,\\ \ell(I) > \ell(J),\\ I\ni y}} |I|
\le |J|\sum_{\nu> 0}\rho^\nu
\le c|J|,
$$
where $0<\rho<1$ is from \eqref{r-rho}.

If $0<\tau< 1$ we apply the same arguments in the estimate of $S_2$
but using the concavity of the function $y^\tau$ instead of H\"{o}lder's inequality.
We have
\begin{align*}
S_2
&\le c\sum_{I\in\II}
\Big(\sum_{j<\ell(I)}\sum_{\substack{\Q\in\QQ_j,\\ V_\Q\cap\Omega_I\ne\emptyset}}
(|I|/|\Q|)|\cc_\Q|\Big)^\tau
\le c\sum_{I\in\II}
\sum_{j<\ell(I)}\sum_{\substack{\Q\in\QQ_j,\\ V_\Q\cap\Omega_I\ne\emptyset}}
(|I|/|\Q|)^\tau|\cc_\Q|^\tau
\\
& \le c\sum_{\Q\in\QQ} |\cc_\Q|^\tau
\sum_{\substack{I\in\II,\\ \ell(I) > \ell(\Q),\\ \Omega_I\cap V_\Q\ne\emptyset}} (|I|/|\Q|)^\tau
\le c\sum_{\Q\in\QQ} |\cc_\Q|^\tau
\le c\|f\|_{\BB^{\alpha, k}_\tau(\FF)}^\tau.
\end{align*}

The above estimates for $S_1$ and $S_2$ and \eqref{B-omega-omega}
yield \eqref{B-E-phi}.

\smallskip

(b) We now show that if $f\in \BB^{\alpha, k}_\tau(E, q)$, $1\le q<\infty$,
then $f\in \BB^{\alpha, k}_\tau(\FF)$ and
\begin{equation}\label{B-phi-Q-E}
\|f\|_{\BB^{\alpha, k}_\tau(\FF)} \le c\|f\|_{\BB^{\alpha, k}_\tau(Q,q)}
\le c\|f\|_{\BB^{\alpha, k}_\tau(E, q)}.
\end{equation}
Indeed, assume that $f\in \BB^{\alpha, k}_\tau(E, q)$.
Then by Theorem~\ref{thm:rep-B} there exists a polynomial $P\in\Pi_k$
such that with unconditional convergence in $\BMO$
\begin{equation*}
f-P = \sum_{Q\in\QQ} b_{Q, q}(f)\varphi_\Q,
\end{equation*}
and by \eqref{NQ=b-norm1} and \eqref{BE-embed-BMO}
\begin{equation*}
\|f\|_{\BB^{\alpha, k}_\tau(Q,q)}
\sim \Big(\sum_{\Q \in \QQ} | b_{\Q,q}(f)|^\tau \Big)^{1/\tau}
\le c\|f\|_{\BB^{\alpha, k}_\tau(E, q)}.
\end{equation*}
Now, using Definition~\ref{def:Besov-B-spline} we conclude that $f\in \BB^{\alpha, k}_\tau(\FF)$
and \eqref{B-phi-Q-E} is valid.

\smallskip

(c) Parts (a), (b) of the proof establish the theorem for any fixed $1\le q<\infty$.
Taking into account that the space $\BB^{\alpha, k}_\tau(\FF)$ does not depend on $q$ this completes the proof.
$\hfill \square$

\begin{rem}\label{rem:Besov}
Observe first that the space $\BB^{\alpha, k}_\tau$ when $\tau\ge 1$
is the usual homogeneous Besov space $\BB^{\alpha}_{\tau\tau}$.
For simplicity of notation we suppress the second index~$\tau$.
Theorem~\ref{thm:B-equiv-norms} shows that
$\|\cdot\|_{\BB^{\alpha, k}_\tau(Q,q)}$ and $\|\cdot\|_{\BB^{\alpha, k}_\tau(E, q)}$ with $1\le q<\infty$
are other equivalent norms in the Besov space $\BB^{\alpha, k}_\tau$, $\tau\ge 1$.
These norms will be more useful for our purposes of spline approximation
than the norm from $\eqref{def-B-norm-modulus}$.

An important difference occurs when $\tau<1$.
The norm $\|\cdot\|_{\BB^{\alpha, k}_\tau}$ from \eqref{def-B-norm-modulus} that involves
the modulus of smoothness $\omega_k(f, t)_\tau$ is hardly usable,
while the norms $\|\cdot\|_{\BB^{\alpha, k}_\tau(Q,q)}$ and $\|\cdot\|_{\BB^{\alpha, k}_\tau(E, q)}$ with $1\le q<\infty$
work very well.
We could have defined the Besov space $\BB^{\alpha, k}_\tau$ for all $\tau>0$ from the outset using the norm
$\|\cdot\|_{\BB^{\alpha, k}_\tau(Q,q)}$.

As is seen above we have introduced $k$ as a parameter and consider $\alpha$ and $k$ independent.
The reason for this is that by allowing $\tau <1$ the spaces $\BB^{\alpha, k}_\tau(Q,q)$ are
nontrivial and different for all $\alpha>0$ and $k\ge 2$.
\end{rem}

According to \thmref{thm:B-equiv-norms} all spaces
$\BB^{\alpha, k}_\tau(Q,q)$, $\BB^{\alpha, k}_\tau(E, q)$, $1\le q<\infty$,
and $\BB^{\alpha, k}_\tau(\Phi)$ (also $\BB^{\alpha, k}_\tau$ if $\tau\ge 1$, see \eqref{B-tau}) are the same;
from now on we shall use the notation $\BB^{\alpha, k}_\tau$ for this Besov space.

\section{Nonlinear spline approximation in $\BMO$}\label{sec:spline-approx}

Assume that $\II$ is a regular multilevel partition of $\RR$.
We denote by $\Phi(\II)$ the collection of all $k-1$ degree B-splines $\varphi_\Q$ generated by $\II$
(see \S \ref{subsec:splines}).
Notice that \emph{$\Phi(\II)$ is not a basis; $\Phi(\II)$ is redundant}.
We consider the nonlinear $n$-term approximation in $\BMO$ from $\Phi(\II)$.

Denote by $\Sigma_n(\II)$ the set of all spline functions $g$ of the form
$$
g = \sum_{\Q \in \Lambda_n} a_\Q \varphi_\Q,
$$
where $a_\Q\in\CC$, $\Lambda_n \subset \QQ(\II)$, $\#\Lambda_n \le n$,
and $\Lambda_n$ may vary with $g$.
We denote by $ \sigma_n(f, \II)_{\BMO}$ the error of $\BMO$-approximation to
$f \in \VMO$ from $\Sigma_n(\II)$:
$$
\sigma_n(f)_{\BMO} = \sigma_n(f, \II)_{\BMO} := \inf_{g \in \Sigma_n(\II)} \|f - g\|_{\BMO}.
$$
Throughout this section we assume as before that
\begin{equation*}
\alpha > 0, \;\; 1/\tau:=\alpha \quad\hbox{and}\quad k\ge 2,
\end{equation*}
and denote by $\BB^{\alpha, k}_\tau$ the Besov space introduced in Section~\ref{sec:Besov-spaces}.

\smallskip
\noindent
{\bf Convention.}
Clearly, if $f\in \BB_\tau^{\alpha, k}$,
then $\|f+P\|_{\BB_\tau^{\alpha, k}}=\|f\|_{\BB_\tau^{\alpha, k}}$ for all $P\in\Pi_k$.
Hence, $\BB_\tau^{\alpha, k}$ consists of equivalence classes modulo $\Pi_k$.
In light of Theorem~\ref{thm:rep-B} for any function $f\in \BB_\tau^{\alpha, k}$
there exists a polynomial $P\in\Pi_k$ such that $f-P\in \VMO$.
From now on we shall assume that each $f\in \BB_\tau^{\alpha, k}$
is the {\em canonical representative} $f-P\in \VMO$ of the equivalence class modulo $\Pi_k$ generated by $f$.
As before (see \S\ref{subsec:BMO}) we identify each $f\in \VMO$ with $f+a$, $a$ is an arbitrary constant.

The following pair of companion Jackson and Bernstein estimates are our main results in this article.
\begin{thm}\label{thm:Jackson} {\bf [Jackson estimate]}
If $f \in \BB^{\alpha, k}_\tau$, then $f\in \VMO$ and
\begin{equation}\label{jackson}
\sigma_n(f, \II)_{\BMO} \le c n^{-\alpha} \| f \|_{\BB^{\alpha, k}_\tau}, \quad n\in\NN,
\end{equation}
with $c>0$ depending only on $\alpha$, $k$, and the parameters of $\II$.
\end{thm}

\begin{thm}\label{thm:bernstein} {\bf [Bernstein estimate]}
If $g \in \Sigma_n(\II)$, $n\in\NN$, then
\begin{equation}\label{bernstein}
\|g\|_{\BB^{\alpha, k}_\tau} \le c n^{\alpha} \|g\|_{\BMO}
\end{equation}
with $c>0$ depending only on $\alpha$, $k$, and the parameters of $\II$.
\end{thm}

\begin{rem}\label{rem:bernstein}
As will be seen from the proof of the Bernstein estimate \eqref{bernstein} $($see \S\ref{sec:bernstein}$)$
the assumption that $g$ being in $\Sigma_n(\II)$ belongs to $C^{k-2}$ is not important;
it is only used that $g\in C$.
Therefore, estimate \eqref{bernstein} is valid for $B$-splines that are only continuous.
\end{rem}

As is well known the companion Jackson and Bernstein estimates \eqref{jackson}, \eqref{bernstein}
imply complete characterization of the approximation spaces associated with spline approximation in $\BMO$.
We next describe this result.

Denote by $K(f, t)$ the $K$-functional associated with $\VMO$ and $\BB^{\alpha, k}_\tau$,
defined for $f\in \VMO$ and $t>0$ by (see e.g. \cite[Chapter 6]{DL})
\begin{equation*}
K(f, t)=K(f, t; \VMO, \BB^{\alpha, k}_\tau)
:= \inf_{g\in \BB^{\alpha, k}_\tau} \big(\|f-g\|_{\BMO} + t\|g\|_{\BB^{\alpha, k}_\tau}\big).
\end{equation*}

The Jackson and Bernstein estimates \eqref{jackson}, \eqref{bernstein}
imply the following direct and inverse estimates:
If $f\in \VMO$, then
\begin{equation}\label{direct-est}
\sigma_n(f)_{\BMO} \le cK\big(f, n^{-\alpha}\big), \quad n\ge 1,
\end{equation}
and
\begin{equation}\label{inver-est}
K(f, n^{-\alpha})
\le cn^{-\alpha}\Big[\Big(\sum_{\nu=1}^{n}\frac{1}{\nu}(\nu^{\alpha}\sigma_\nu(f,\II)_{\BMO})^\mu\Big)^{1/\mu} +\|f\|_{\BMO}\Big],
\quad n\ge 1,
\end{equation}
where $\mu:=\min\{\tau, 1\}$.

The proofs of \eqref{direct-est}, \eqref{inver-est} are standard,
see e.g. \cite[Chapter 7, Theorem 5.1]{DL}.

We define the approximation space $A^\gamma_q(\BMO, \II)$ generated by
nonlinear $n$-term approximation in $\BMO$ from B-spliens
to be the set of all functions $f\in BMO$ such that
\begin{equation*}
\|f\|_{A^\gamma_q(\BMO, \II)}
:= \|f\|_{\BMO} + \Big(\sum_{n=1}^\infty \big(n^\gamma \sigma_n(f, \II)_{\BMO}\big)^q\frac{1}{n}\Big)^{1/q} <\infty,
\end{equation*}
with the usual modification when $q=\infty$.

The following characterization of the approximation spaces $A^\gamma_q(\BMO, \II)$
is immediate from inequalities \eqref{direct-est}, \eqref{inver-est}:

\begin{thm}\label{thm:app-sp}
If $0<\gamma<\alpha$ and $0<q\le \infty$, then
\begin{equation*}
A^\gamma_q(\BMO, \II) = \big(\VMO, \BB^{\alpha, k}_\tau\big)_{\frac{\gamma}{\alpha}, q}.
\end{equation*}
In particular, if $f\in \VMO$, then
\begin{equation*}
K(f, t^\alpha) = O(t^\gamma)
\quad\hbox{if and only if}\quad
\sigma_n(f)_{\BMO}=O(n^{-\gamma}).
\end{equation*}
Above $(X_0, X_1)_{\lambda, q}$ stands for the real interpolation space
induced by two spaces $X_0$, $X_1$,
see e.g. \cite[Chapter 6, \S7]{DL}.
\end{thm}

\section{\!\!Comparison with spline approximation in other spaces}\label{sec:comparison}

\subsection{\!Comparison between spline approximation in \boldmath $\BMO$ and $L^\infty$}\label{subsec:comparison}

Here we clarify the differences and similarities between
nonlinear $n$-term approximation from B-splines in $\BMO$ and in the uniform norm.

Denote by $ \sigma_n(f, \II)_\infty$ the error of $L^\infty$-approximation to
$f$ from $\Sigma_n(\II)$:
$$
\sigma_n(f, \II)_\infty := \inf_{S \in \Sigma_n(\II)} \|f - S\|_\infty.
$$
We denote by $\BB_{\tau}^{\alpha, k}$ the Besov space introduced in Section~\ref{sec:Besov-spaces}.
The following theorem follows from the Jackson estimates in \cite[Theorem~4.1]{DahP} or \cite[Theorem~4.1]{KPS}.

\begin{thm}\label{thm:Jackson-C} {\bf [Jackson estimate]}
If $f \in \BB_{\tau}^{\alpha, k}$, $\alpha\ge 1$, $1/\tau=\alpha$, $k\ge 2$,
then $f$ is continuous on $\RR$, $\lim_{|x|\to \infty} f(x)=0$, and
\begin{equation}\label{jackson-C}
\sigma_n(f, \II)_\infty \le c n^{-\alpha} \| f \|_{\BB^{\alpha, k}_\tau}, \quad n\in\NN,
\end{equation}
with $c$ depending only on $\alpha$ and the parameters of $\II$.
\end{thm}

Observe that the B-spaces used for nonlinear spline approximation in $L^\infty$ in \cite{DahP, KPS}
are somewhat different because the approximation there takes place in dimension $d\ge 1$ or $d=2$.
However, from (3.5) in \cite{DahP} or (2.15) in \cite{KPS} it is clear that these spaces
are the same as $\BB_{\tau}^{\alpha, k}$ in dimension $d=1$.

\begin{thm}\label{thm:bernstein-C} {\bf [Bernstein estimate]}
Let $\alpha >0$ and $1/\tau=\alpha$.
If $S \in \Sigma_n(\II)$, $n\in\NN$, then
\begin{equation}\label{bernstein-C}
\| S \|_{\BB_\tau^{\alpha, k}} \le c n^{\alpha} \| S \|_\infty,
\end{equation}
with $c$ depending only on $\alpha$ and the parameters of $\II$.
\end{thm}

This theorem follows by the Bernstein estimates in \cite[Theorem~4.2]{DahP}, see also \cite[Theorem~4.2]{KPS}.

Several clarifying remarks are in order:

(1) The Besov space $\BB^{\alpha, k}_\tau$ is obviously not embedded in $L^\infty$ when $\alpha<1$.
For this reason the Jackson estimate \eqref{jackson-C} is not valid when $\alpha<1$.
At the same time the Jackson estimate \eqref{jackson} holds for all $\alpha>0$.

(2) Since $\|S\|_{\BMO} \le \|S\|_\infty$ the Bernstein estimate \eqref{bernstein-C}
is a consequence of the Bernstein estimate \eqref{bernstein}.
Similarly, the Jackson estimate \eqref{jackson} follows by \eqref{jackson-C}
in the case when $\alpha\ge 1$.

(3) It is interesting that in the case when $\alpha\ge 1$ the same Besov spaces $\BB^{\alpha, k}_\tau$
work for spline approximation in both $\BMO$ and $L^\infty$.

(4) Algorithms for nonlinear $n$-term approximation from linear B-splines are developed in \cite{DP}
and in more general settings in \cite{DahP}.
The results in \cite{DahP, KPS} use ideas that originate in earlier development of
spline approximation in $L^\infty$ in \cite{DPY}.

\subsection{Comparison between spline approximation in \boldmath $\BMO$ and $L^p$, $p<\infty$}\label{subsec:comparison-2}

There is a principle difference between nonlinear spline approximation in $\BMO(\RR)$ (or $L^\infty(\RR)$)
and in $L^p(\RR)$, $1\le p<\infty$, that we would like to clarify here.

To be specific, denote by $S(k, n)$ the set of all piecewise polynomials functions on $\RR$
of degree $k-1$ with $n+1$ free knots,
that is, $S\in S(k, n)$ if there exist points
$-\infty<x_0<x_1<\cdots<x_n<\infty$
and polynomials $P_j\in \Pi_k$, $j=1, \dots, n$, such that
\begin{equation}\label{S-k-n}
S=\sum_{j=1}^n \ONE_{I_j}\cdot P_j,
\quad I_j:=[x_{j-1}, x_j).
\end{equation}
Here the knots $\{x_j\}$ are allowed to vary with $S$.
Hence $S(k, n)$ is nonlinear. No multilevel partition is assumed.
Given $f\in L^p(\RR)$, $1\le p<\infty$, define
\begin{equation*}
S_n^k(f)_p:= \inf_{S\in S(k, n)}\|f-S\|_{p}.
\end{equation*}
One is interested in characterizing the approximation spaces associate to this approximation process.
The Besov spaces
\begin{equation*}
\BB^{\alpha, k}_\tau, \quad \alpha>0,\;\; 1/\tau:=\alpha+1/p,\;\; k\ge 1,
\end{equation*}
naturally appear in this sort of problems.
These spaces are standardly defined \cite{P} by the following norm using muduli of smoothness as in \eqref{def-B-norm-modulus}
for $0<\tau<\infty$:
\begin{equation*}
\|f\|_{\BB_\tau^{\alpha, k}} := \Big(\int_0^\infty [t^{-\alpha}\omega_k(f, t)_\tau]^\tau\frac{dt}{t}\Big)^{1/\tau}.
\end{equation*}
However, the definition of $\BB^{\alpha, k}_\tau$ can be modified as in Definition~\ref{def:B-2}
when $\tau<1$, $q<p$, see \cite{DP, KP}.

The following Jackson and Bernstein estimates have been established in \cite{P}:
If $f\in \BB_\tau^{\alpha, k}$, then $f\in L^p$ and
\begin{equation*}
S_n^k(f)_p \le cn^{-\alpha}\|f\|_{\BB_\tau^{\alpha, k}}
\end{equation*}
and for any $S\in S(k, n)$
\begin{equation*}
\|S\|_{\BB_\tau^{\alpha, k}} \le cn^\alpha\|S\|_p.
\end{equation*}
Clearly, these two estimates allow to completely characterize the associated to
$\{S_n^k(f)_p\}$ approximation spaces.

\smallskip
\noindent
{\bf Discussion.}
As the following remarks show the nature of nonlinear spline (piecewise polynomials) approximation
in $\BMO(\RR)$ (or $L^\infty(\RR)$) and in $L^p(\RR)$, $p<\infty$, is totaly different.

When approximating from piecewise polynomials in $L^p$, $1\le p<\infty$,
there is no need to assume any smoothness or continuity,
we simply work with discontinuous piecewise polynomials, see \eqref{S-k-n}.
The point is that smooth piecewise polynomials and discontinuous piecewise polynomials produce the same rates of approximation for $p<\infty$.
More importantly, if $S_I:= \ONE_I\cdot P$ for some compact interval $I$
and a polynomial $P\in\Pi_k$, $P\not\equiv 0$, then
\begin{equation*}
\|S_I\|_{\BB_\tau^{\alpha, k}} <\infty, \quad \forall \alpha>0, \;\; 1/\tau=\alpha+1/p.
\end{equation*}
In other words this piecewise polynomial function is infinitely smooth
in the scale of the Besov spaces $\BB_\tau^{\alpha, k}$.

In contrast, it is easy to see for $S_I\notin C(\RR)$ that
\begin{equation*}
\|S_I\|_{\BB_\tau^{\alpha, k}} =\infty
\quad\hbox{for any $\alpha>0$ whenever $1/\tau=\alpha$, i.e. $p=\infty$.}
\end{equation*}
Furthermore, the Bernstein estimate \eqref{bernstein} cannot be true
for piecewise polynomials $S$ of the form \eqref{S-k-n} even if $S\in C^{k-2}$.
We next clarify this claim with the following simple example.
Let $S(x)=1$ for $x\in [0, 1]$, $S(x)=0$ for $x\in (-\infty, -\varepsilon]\cup [1+\varepsilon, \infty)$,
and $S$ is linear and continuous on $[-\varepsilon, 0]$ and $[1, 1+\varepsilon]$,
where $\varepsilon>0$ is sufficiently small.
It is readily seen that
$\omega_2(S, t)_\tau^\tau \sim t$ for $\varepsilon \le t\le 1/2$
and hence
\begin{equation*}
\|S\|_{\BB^{\alpha, 2}_\tau}^\tau \ge c\int_\varepsilon^{1/2}\frac{dt}{t} \ge c \ln(1/\varepsilon).
\end{equation*}
Therefore, the Bernstein estimate \eqref{bernstein} is not valid for piecewise polynomials
of the form \eqref{S-k-n} even if they are smooth.
This is the reason for considering nonlinear approximation from splines generated
by a hierarchy of B-splines.

\section{Proof of Theorem~\ref{thm:Jackson}}\label{sec:jackson}

In this section we prove the Jackson estimate \eqref{jackson}.
We shall derive this estimate from an estimate for approximation in general sequence spaces.

\subsection{Jackson inequality for nonlinear approximations in $\ff^{0 q}_\infty$}

Here we consider nonlinear $n$-term approximation from finitely supported sequences
in the spaces $\ff^{0 q}_\infty$ in a general setting developed in \cite[\S7.2]{IP1}.

\begin{defn}\label{nested_structure}
Let $\cX=\cup_{m=-\infty}^\infty\cX_m$ be a countable multilevel index set.
With every $\xi\in\cX$ we associate an open set $U_\xi\subset\RR$.
We call $\{U_\xi : \xi\in\cX\}$ a \emph{nested structure associate with $\cX$} if there exists constant $\lambda\ge 1$ such that:
\begin{enumerate}
\renewcommand{\theenumi}{\alph{enumi}} 
\renewcommand{\labelenumi}{$(\theenumi)$}
\item $|\RR\backslash\bigcup_{\xi\in\cX_m}U_\xi|=0\quad\forall m\in\ZZ$;
\item If $\eta\in\cX_m, \xi\in\cX_\nu$ and $m\ge \nu$ then either $U_\eta\subset U_\xi$ or $U_\eta\cap U_\xi=\emptyset$;
\item For every $\xi\in\cX_\nu, m< \nu,$ there is unique $\eta\in\cX_m$ such that $U_\xi\subset U_\eta$;
\item $|U_\eta|\le \lambda|U_\zeta| \quad\forall\eta,\zeta\in\cX_m, \forall m\in\ZZ$;
\item Every $U_\xi$, $\xi\in\cX_m$, has at least two children, i.e. there are $\eta,\zeta\in\cX_{m+1}$
      such that $U_\eta\subset U_\xi$, $U_\zeta\subset U_\xi$, $\eta\ne\zeta$.
\end{enumerate}
We also assume that there exist one $\fD^1$ or two $\fD^1,\fD^2$ disjoint subsets of $\RR$
with the properties:
$(i)$ $|\RR\setminus\cup_{j=1}^K\fD^j|=0$;
$(ii)$ for every $\xi\in\cX$ there is $j,~1~\le~j~\le~K$, such that $U_\xi\subset \fD^j$;
$(iii)$ if the sets $U_\eta,U_\zeta$ are contained in $\fD^j$ for some $1~\le~j~\le~K$,
then there exists $U_\xi\subset\fD^j$ such that $U_\eta\subset U_\xi$, $U_\zeta\subset U_\xi$.
Here $K=1$ or $K=2$.
\end{defn}

\begin{rem}
Conditions $(a)$--$(c)$ imply that for every $\xi\in\cX_m$ we have
\begin{equation*}
|U_\xi|=\sum_{\substack{\eta\in\cX_{m+j}\\ U_\eta\subset U_\xi}}|U_\eta|,
\quad \forall j\in\NN.
\end{equation*}

Conditions $(d)$--$(e)$ imply that there exists $\rho\in(0,1)$, namely $\rho=\lambda/(\lambda+1)$, such that
\begin{equation}\label{dyad_eq:6}
|U_\eta|\le \rho |U_\xi|\quad\forall m\in\ZZ,~\xi\in\cX_m,~\eta\in\cX_{m+1},~U_\eta\subset U_\xi.
\end{equation}
\end{rem}

Recall that by definition the sequence spaces $\ell^\tau=\ell^\tau(\cX)$, $0<\tau\le\infty$,
consists of all sequences $\{h_\xi : \xi\in\cX\}$ such that
\begin{equation*}
\|h\|_{\ell^\tau}
:= \bigg(\sum_{\xi\in \cX}|h_\xi|^\tau\bigg)^{1/\tau}<\infty.
\end{equation*}
We define the sequence spaces $\fg^q=\fg^q(\cX)$, $0<q\le\infty$, as the set of all sequences $\{h_\xi : \xi\in\cX\}$ such that
\begin{equation*}
\|h\|_{\fg^q}
:= \sup_{\xi\in\cX}\bigg(\sum_{U_\eta\subset U_\xi}|h_\eta|^q \frac{|U_\eta|}{|U_\xi|}\bigg)^{1/q}<\infty.
\end{equation*}
Note that the definition of $\ell^\tau(\cX)$ does not need a nested structure associate with $\cX$,
while the definition of $\fg^q(\cX)$, $q<\infty$, requires such kind of structure.
Of course, for $q=\tau=\infty$ we have $\fg^\infty=\ell^\infty$.

The best nonlinear approximation of $h\in\fg^q$ from sequences with at most $n$ non-zero elements is given by
\begin{equation*}
\sigma_n(h)_{\fg^q}:=\inf_{|\supp~\bar{h}|\le n}\|h-\bar{h}\|_{\fg^q}
=\inf_{\substack{\Lambda_n\subset\cX\\ \#\Lambda_n\le n}}\sup_{\xi\in\cX}\Bigg(
\sum_{\substack{U_\eta\subset U_\xi\\ \eta\notin\Lambda_n}}|h_\eta|^q \frac{|U_\eta|}{|U_\xi|}\Bigg)^{1/q}.
\end{equation*}

\begin{thm}[Jackson inequality]\label{thm:Jackson_seq}
Let $0<\tau<\infty$ and $0<q\le\infty$.
Assume $\{U_\xi : \xi\in\cX\}$ is a nested structure associate with $\cX$.
There exists a constant $c=c(\tau,q,\lambda)$
such that for any $h\in\ell^{\tau}$ we have $h\in\fg^q$ and
\begin{equation*}
\sigma_n(h)_{\fg^q}\le c n^{-1/\tau}\|h\|_{\ell^\tau}, \quad n\in\NN.
\end{equation*}
\end{thm}

\subsection{Proof of \thmref{thm:Jackson}}

Let the support $Q\in\QQ_m$ of $\varphi$ be $Q=[x_{m,j},x_{m,j+k}]$.
We use $\Q$ in the place of $\xi$,
$\QQ_m$ in the place of $\cX_m$, $m\in\ZZ$, $\QQ$ in the place of $\cX$ and
the open interval $I_{Q}:=(x_{m,j},x_{m,j+1})$ in the place of $U_\xi$
from Definition~\ref{nested_structure}.
If $\II$ is a regular multilevel partition of $\RR$
then the system of intervals $\{I_{Q} : Q\in\QQ\}$ is a nested structure associate with $\QQ=\cup_{m=-\infty}^\infty \QQ_m$.
It satisfies all conditions of Definition~\ref{nested_structure}.

In order to determine $K$ from Definition~\ref{nested_structure} we consider all knots $\{x_{m,j} : m,j\in\ZZ\}$.
In view of \eqref{dyad_eq:6} it is not possible to have two different points $y_1,y_2\in\RR$ that are common knots for all levels $m\in\ZZ$.
If all levels have one common knot $y_1$ then $K=2$, $\fD^1=(-\infty,y_1)$, $\fD^2=(y_1,\infty)$.
If there is no common knot for all levels then $K=1$, $\fD^1=(-\infty,\infty)$.

\begin{lem}\label{lem:BMO_seq}
If $\{a_Q\}_{Q\in\QQ}\in\fg^1(\QQ)$ and $f=\sum_{Q\in\QQ}a_Q\varphi_Q\in\VMO(\RR)$ then
\begin{equation*}
\|f\|_{\BMO}\le c \|\{a_Q\}\|_{\fg^1}.
\end{equation*}
\end{lem}
\proof
Denote
\begin{equation*}
f_\nu:=\sum_{j>\nu}\sum_{\Q\in\QQ_j} a_\Q\varphi_\Q,
\quad \nu\in\ZZ.
\end{equation*}
We claim that
\begin{equation}\label{est-BMO-fk2}
\|f_\nu\|_{\BMO}
\le c\sup_{\Q\in\QQ_j, j>\nu}\sum_{I_{\Q'}\subset I_\Q}|a_{\Q'}| \frac{|I_{\Q'}|}{|I_\Q|}
=:\|\{a_\Q\}\|_{\fg^1(\nu)}.
\end{equation}
Let $J$ be an arbitrary compact interval in $\RR$.
Then there exist $m\in\ZZ$ such that
if $\Q\in\QQ_m$ and $\Q\cap J\ne \emptyset$, then $|\Q|\sim |J|$  and $\Q\subset 2J$.
Denote $\tilde{J}:=2J$.

We may assume that $\nu<m$.
We split $f_\nu$ into two: $f_\nu=f_m+(f_\nu-f_m)$.
Using \eqref{equiv-norms} we get
\begin{multline}\label{BMO-fm2}
\frac{1}{|J|}\int_J|f_m(x)|dx
\le \frac{1}{|J|}\int_{J}\sum_{\Q'\subset \tilde{J}} |a_{\Q'}\varphi_{\Q'}(x)|dx
\le \frac{c}{|\tilde{J}|}\int_{\tilde{J}}\sum_{\Q'\subset \tilde{J}} |a_{\Q'}\varphi_{\Q'}(x)|dx\\
\le c\sum_{\Q\in\QQ_m, \Q\cap \tilde{J}\ne\emptyset}\sup_{\Q\in\QQ_j, j>m}\sum_{I_{\Q'}\subset I_\Q}|a_{\Q'}| \frac{|I_{\Q'}|}{|I_\Q|}
\le c\|\{a_\Q\}\|_{\fg^1(\nu)}.
\end{multline}

Fix $y\in J$. Denote $F_{\nu m}:=f_\nu-f_m$.
We claim that
\begin{equation}\label{est-Fkm2}
\frac{1}{|J|}\int_J|F_{\nu m}(x)-F_{\nu m}(y)|dx
\le c\|\{a_\Q\}\|_{\fg^1(\nu)}.
\end{equation}
Indeed, let $\Q\in\QQ_j$, $j\le m$,
and assume $\Q\cap J\ne \emptyset$.
Then for $x\in \Q$
\begin{equation*}
|\varphi_\Q(x)- \varphi_\Q(y)| \le |x-y||\varphi_\Q'(\xi)| \le c|x-y||\Q|^{-1},
\quad \xi\in (x, y).
\end{equation*}
For every $x\in J$ using the above we get
\begin{align*}
|F_{\nu m}(x)-F_{\nu m}(y)|
&\le \sum_{j=\nu+1}^m\sum_{\Q\in\QQ_j, \Q\ni x} |c_\Q||\varphi_\Q(x)-\varphi_\Q(y)|
\\
&\le c\|\{a_\Q\}\|_{\fg^1(\nu)}\sum_{j=\nu+1}^m\sum_{\Q\in\QQ_j, \Q\ni x} |\varphi_\Q(x)-\varphi_\Q(y)|
\\
&\le c\|\{a_\Q\}\|_{\fg^1(\nu)} \sum_{j=\nu+1}^m\sum_{\Q\in\QQ_j, \Q\ni x} |J||\Q|^{-1}
\le c\|\{a_\Q\}\|_{\fg^1(\nu)}.
\end{align*}
In the last inequality we used \eqref{r-rho} and that every point $x$ is contained in $k$ different $\Q$'s from every level.
Estimate \eqref{est-Fkm2} follows readily from the above inequalities.

From \eqref{BMO-fm2} and \eqref{est-Fkm2} it follows that
\begin{equation*}
\frac{1}{|J|}\int_J|f_\nu(x)-F_{\nu m}(y)|dx
\le c\|\{a_\Q\}\|_{\fg^1(\nu)},
\end{equation*}
which implies \eqref{est-BMO-fk2}.

Now, condition $f\in\VMO(\RR)$  implies $\lim_{\nu\to -\infty}\|f_\nu-f\|_{\BMO}=0$ and,
hence, \eqref{est-BMO-fk2} proves the lemma.
$\hfill \square$

\smallskip
\noindent
{\it Completion of the proof of \thmref{thm:Jackson}.}
Let $h=\{a_Q\}_{Q\in\QQ}\in\ell^\tau(\QQ)$. Then \thmref{thm:Jackson_seq} with $n=1$ implies $h\in\fg^1$ and \thmref{thm:embed-B-BMO} gives
$f=\sum_{Q\in\QQ}a_Q\varphi_Q\in\VMO$. Using \lemref{lem:BMO_seq} and once more time \thmref{thm:Jackson_seq} we obtain
\begin{equation*}
\sigma_n(f)_{\BMO(\RR)} \le c \sigma_n(h)_{\fg^1(\QQ)}\le c n^{-1/\tau}\|h\|_{\ell^\tau(\QQ)}\le c n^{-\alpha} \| f \|_{\BB^{\alpha, k}_\tau},\quad n\in\NN.
\end{equation*}
This proves \thmref{thm:Jackson}.
$\hfill \square$

\section{Proof of Theorem~\ref{thm:bernstein} }\label{sec:bernstein}

For the proof of the Bernstein estimate \eqref{bernstein}
we shall need Lemma~\ref{lem:poly-norms} and the following

\begin{lem}\label{lem:L1-BMO}
Let $f\in \BMO$ and $f(x)=0$ for $x\in J\setminus I$,
where $I, J$ are two intervals so that $I\subset J$ and $|J|=(1+\delta)|I|$, $\delta>0$.
If $1\le \tau <\infty$, then
\begin{equation}\label{Lp-BMO}
\int_I|f(x)|^\tau dx
\le (1+\delta^{-1})^\tau\int_J|f(x)-\avg_Jf|^\tau dx
\le c|J|\|f\|_{\BMO}^\tau,
\end{equation}
where $c=c(\delta, \tau)$.
\end{lem}
\proof
Using the hypothesis of the lemma we have
\begin{align*}
&\Big(\frac{1}{|I|}\int_I|f(x)|^\tau dx\Big)^{1/\tau}
\le \Big(\frac{1}{|I|}\int_I|f(x)-\avg_Jf|^\tau dx\Big)^{1/\tau} + |\avg_Jf|\\
&\le \Big(\frac{1}{|I|}\int_I|f(x)-\avg_Jf|^\tau dx\Big)^{1/\tau} + \frac{1}{|J|}\int_J|f(x)|dx\\
&\le \Big(\frac{1}{|I|}\int_J|f(x)-\avg_Jf|^\tau dx\Big)^{1/\tau} + \frac{|I|}{|J|}\frac{1}{|I|}\int_I|f(x)|dx\\
&\le \Big(\frac{1}{|I|}\int_J|f(x)-\avg_Jf|^\tau dx\Big)^{1/\tau} + \frac{1}{1+\delta}\Big(\frac{1}{|I|}\int_I|f(x)|^\tau dx\Big)^{1/\tau}
\end{align*}
with H\"{o}lder's inequality applied in the last inequality. This proves the first inequality in \eqref{Lp-BMO},
while the second inequality follows from \eqref{equiv-norms-BMO}.
$\hfill \square$

\medskip
\noindent
{\it Proof of Theorem~\ref{thm:bernstein}.}
Given a 
partition $\II$ we
set $t_m=\lambda\sup_{I\in\II_m}|\Omega_I|$, $m\in\ZZ$,  with $\lambda$ from \eqref{cond-rho}.
From \eqref{cond-rho} we infer that $t_m\lambda^{-2}\le|\Omega_I|\le t_m\lambda^{-1}$ for every $I\in\II_m$
and from \eqref{r-rho} we obtain $r t_m\le t_{m+1}\le\rho t_m$, $m\in\ZZ$.

Let $g=\sum_{\Q\in\Lambda_n}\cc_\Q\varphi_\Q\in \Sigma_n$.
Denote by $x_0<x_1<\cdots<x_N$, $N \le (k+1)n$, the knots of $\varphi_\Q$, $\Q\in\Lambda_n$, in increasing order
(if two B-splines have a common knot then it appears only once in this sequence).
There exist polynomials $P_j\in \Pi_{k-1}$, $j=1,\dots,N$, such that
\begin{equation}\label{rep-g}
g=\sum_{j=1}^N \ONE_{J_j}\cdot P_j,
\quad J_j:=[x_{j-1}, x_j],
\quad g\in C(\RR).
\end{equation}
In fact $g\in C^{k-2}(\RR)$.
Also, denote $J_0:=(-\infty, x_0]$ and $J_{N+1}:=[x_N, \infty)$.
Note that the points $x_0<x_1<\cdots<x_N$ are among the knots of partition $\II$
and two consecutive $x_{j-1}, x_j$  may belong to different levels of $\II$.

(a) Let $1\le\tau<\infty$.
 For every $m\in\ZZ$ and every $I\in\II_m$ we shall establish the estimate
\begin{equation}\label{est_E_I}
E_k(g,\Omega_I)_\tau^\tau\le c\bigg(\sum_{J_\nu\subset\Omega_I} |J_\nu|
+\sum_{\substack{0<|J_\mu\cap\Omega_I|<|\Omega_I|\\J_\mu\setminus\Omega_I\ne\emptyset}}
\min\Big\{|J_\mu|,\frac{t_m^{1+\tau}}{|J_\mu|^\tau}\Big\}\bigg)\|g\|_{\BMO}^\tau.
\end{equation}
The second sum in \eqref{est_E_I}, where the summation is on $J_\mu$,  contains 0, 1 or 2 terms.
Every such term represents an interval $J_\mu$ which \emph{partially} covers $\Omega_I$
and contains \emph{in its interior} one of the end points of $\Omega_I$.
Note that $J_0$ and $J_{N+1}$ have length $\infty$ and if they are among $J_\mu$'s then the corresponding term is $0$.
We shall term the intervals $J_j, j=1,\dots,N,$ with $|J_j|>t_m$ \emph{big intervals} for the level $m$.

For the proof of \eqref{est_E_I} we consider five cases depending on the position
of the interval $\Omega_I$, $I\in\II_m$, relative to the intervals $J_j$, $j=0,1,\dots,N+1$, of $g$.

{\em Case 1. $\Omega_I$ is a subset to one of the intervals $J_j$, $j=0,1,\dots,N+1$.}

In this case \eqref{est_E_I} is trivially satisfied because $E_k(g,\Omega_I)_\tau=0$.

{\em Case 2. Both end points of $\Omega_I$ are either among the knots $x_j, j=1,\dots,N$, or are in the interior of intervals $J_\mu$ with $|J_\mu|\le t_m$.}

Note that $\min\{|J_\mu|,t_m^{1+\tau}/|J_\mu|^\tau\}=|J_\mu|$ if $J_\mu$ belongs to the second sum in \eqref{est_E_I}.
Set $J^\star=\cup_{|J_\nu\cap\Omega_I|>0}J_\nu$. Then \eqref{est_E_I} follows from
\begin{equation*}
E_k(g,\Omega_I)_\tau^\tau\le E_1(g,J^\star)_\tau^\tau\le |J^\star|\|g\|_{\BMO}^\tau=\sum_{|J_\nu\cap\Omega_I|>0}|J_\nu|\|g\|_{\BMO}^\tau.
\end{equation*}

{\em Case 3. Both end points of $\Omega_I$ are in the interior of two intervals $J_\mu$ with $|J_\mu|> t_m$
and there is only one knot among $x_j, j=1,\dots,N$, in the interior of $\Omega_I$.}

Let $x_j$ be the knot of $g$ in the interior of $\Omega_I$. Then the two big intervals covering the ends of $\Omega_I$ are $J_j$ and $J_{j+1}$. We have
\begin{align*}
E_k(g,\Omega_I)_\tau^\tau\le E_1(g,\Omega_I)_\tau^\tau &\le c|\Omega_I|\big(|\Omega_I|\|g'\|_{L^\infty(J_j\cup J_{j+1})}\big)^\tau
\\
&\le c|\Omega_I|^{1+\tau}\big(\|g'\|_{L^\infty(J_j)}^\tau+ \|g'\|_{L^\infty(J_{j+1})}^\tau\big).
\end{align*}
Further, using Lemma~\ref{lem:poly-norms} we get for $J=J_j$ or $J=J_{j+1}$
\begin{align}\label{est-BMO1}
\|g'\|_{L^\infty(J)}
&= \|(g-\avg_{J}g)'\|_{L^\infty(J)}\nonumber
\\
&\le \frac{c}{|J|}\|g-\avg_{J}g\|_{L^\infty(J)}
\le \frac{c}{|J|^2}\|g-\avg_{J}g\|_{L^1(J)}
\\
&=\frac{c}{|J|}\cdot\frac{1}{|J|}\int_{J}|g(x)-\avg_{J}g|dx
\le \frac{c}{|J|}\|g\|_{\BMO}.\nonumber
\end{align}
Therefore
\begin{equation}\label{est-Y1}
E_k(g,\Omega_I)_\tau^\tau
\le ct_m^{1+\tau}\sum_{\substack{0<|J_\mu\cap\Omega_I|<|\Omega_I|\\J_\mu\setminus\Omega_I\ne\emptyset}} \frac{1}{|J_\mu|^\tau}\|g\|_{\BMO}^\tau.
\end{equation}
Note that \eqref{est-Y1} also holds if one of $J_j$ and $J_{j+1}$ is unbounded because $g=0$ on this interval.
Inequality \eqref{est_E_I} reduces to \eqref{est-Y1} in this case.

{\em Case 4. Both end points of $\Omega_I$ are in the interior of two intervals $J_\mu$ with $|J_\mu|> t_m$
and there are at least two knots among $x_j, j=1,\dots,N$, inside $\Omega_I$.}

Let $x_{j_1},\dots,x_{j_2-1}$, $x_{j_1}<x_{j_2-1}$, be the knots of $g$ in the interior of $\Omega_I$.
Then the intervals $J_\mu$ from the second sum in \eqref{est_E_I} are $J_{j_1}$ are $J_{j_2}$.
Denote $J=[x_{j_1},x_{j_2-1}], J^\star=J_{j_1}\cup J\cup J_{j_2}$, and
\begin{equation*}
V=\sum_{\Q\in\Lambda_n,~\#V_\Q\cap J\ge 2}\cc_\Q\varphi_\Q.
\end{equation*}
Recall $V_\Q$ denotes the knots of $\varphi_\Q$.
Note that $|J|<|\Omega_I|$ and $|J_1|,|J_2|>\lambda |\Omega_I|$.
Hence $\#V_\Q\cap J\ge 2$ implies $\#V_\Q\cap J=k+1$ (for $\Q\in\Lambda_n$),
i.e. if a B-spline is involved in the sum of $V$, then it has all knots in $J$.
The same argument shows that if $g$ has two B-splines from \emph{different} levels with only one knot in $J$,
then this is a common knot.
(Assuming the contrary, from the refinement property of $\II$ we get that
the higher level contains an interval of length $<|\Omega_I|$ and hence it cannot contain an interval of length $>\lambda |\Omega_I|$.
This contradicts the assumption that a B-spline from this level has a single knot in $J$.)
Denote by $y$ the common knot of all B-splines of $g$ with only one knot in the interior of $J$ (if there are such splines).
Hence we can write
\begin{equation*}
g(x)=U(x)+V(x),\quad x\in J^\star,
\end{equation*}
where
\begin{equation*}
U=\ONE_{[x_{j_1-1},y]}P_{j_1}+\ONE_{[y,x_{j_2}]}P_{j_2}, \quad P_{j_1}, P_{j_2}\in \Pi_{k-1}, \quad U\in C(J).
\end{equation*}
Here $P_{j_1}, P_{j_2}$ are from \eqref{rep-g}.
In case there are no B-splines of $g$ with only one knot in the interior of $J$
the above representation holds with $P_{j_1}=P_{j_2}$ and $y$ being any of $x_{j_1},\dots,x_{j_2-1}$.

We have
\begin{equation}\label{est-g}
E_k(g,\Omega_I)_\tau^\tau\le cE_k(U,\Omega_I)_\tau^\tau+cE_k(V,\Omega_I)_\tau^\tau.
\end{equation}

To estimate the best approximation of $V$ we use Lemma~\ref{lem:L1-BMO} and $3J\subset J^\star$ to obtain
\begin{align}\label{est-V}
E_k(V,\Omega_I)_\tau^\tau&\le c\|V\|_\tau^\tau
=c\int_{J}|V(x)|^\tau dx \le c\int_{3J}|V(x)-\avg_{3J}V|^\tau dx \notag \\
&\le c\int_{3J}|g(x)-\avg_{3J}g|^\tau dx + c\int_{3J}|U(x)-\avg_{3J}U|^\tau dx \\
&\le c|J|\|g\|_{\BMO}^\tau + c|J|\big(|J|\|U'\|_{L^\infty(J^\star)}\big)^\tau. \notag
\end{align}
To estimate the best approximation of $U$ we write
\begin{equation}\label{est-U}
E_k(U,\Omega_I)_\tau^\tau\le E_1(U,\Omega_I)_\tau^\tau \le c|\Omega_I|\big(|\Omega_I|\|U'\|_{L^\infty(J^\star)}\big)^\tau
\end{equation}
Inserting \eqref{est-V} and \eqref{est-U} in \eqref{est-g} we obtain
\begin{equation}\label{est-g2}
E_k(g,\Omega_I)_\tau^\tau\le c\sum_{J_\nu\subset\Omega_I} |J_\nu|\|g\|_{\BMO}^\tau + ct_m^{1+\tau}\|U'\|_{L^\infty(J^\star)}^\tau.
\end{equation}

For the estimate of $U'$ we use Lemma~\ref{lem:poly-norms} and \eqref{est-BMO1} with $J=J_{j_1}$ to obtain
\begin{equation}\label{est-Up}
\|U'\|_{L^\infty([x_{j_1-1},y])}\le c\|U'\|_{L^\infty(J_{j_1})}=c\|g'\|_{L^\infty(J_{j_1})}\le \frac{c}{|J_{j_1}|}\|g\|_{\BMO}.
\end{equation}
Similarly for the interval $[y,x_{j_2}]\supset J_{j_2}$.

Combining \eqref{est-g2} and \eqref{est-Up} we obtain \eqref{est_E_I}.
Estimate \eqref{est_E_I} is also valid in the case when one or both intervals $J_{j_1}$, $J_{j_2}$ are unbounded because $g=U=0$ here.

{\em Case 5.  One of the end points of $\Omega_I$ is among the knots $x_j, j=1,\dots,N$, or is in the interior of an interval $J_\mu$ with $|J_\mu|\le t_m$
and the other end point is in the interior of an interval $J_\mu$ with $|J_\mu|> t_m$.}

This is a  simplified version of Case 4.
Let $x_{j_1},\dots,x_{j_2-1}$, $x_{j_1}\le x_{j_2-1}$, be the knots of $g$ in the interior of $\Omega_I$.
Without loss of generality, let $|J_{j_1}|> t_m$ be the big interval containing the left end of $I$.
Then for the right end of $I$ belongs to $J_{j_2}=[x_{j_2-1},x_{j_2}]$ and $|J_{j_2}|\le t_m$.
Set $J=[x_{j_1},x_{j_2}], J^\star=J_{j_1}\cup J$, and
\begin{equation*}
U=\ONE_{J^\star}P_{j_1},\quad P_{j_1}\in\Pi_k;\qquad V(x)=g(x)-U(x),\quad x\in J^\star.
\end{equation*}
Thus, the polynomial $U$ coincides with $g$ on $J_{j_1}$ and the spline $V$ is zero on  $J_{j_1}$. We have
\begin{equation}\label{est-V5}
E_k(g,\Omega_I)_\tau^\tau=E_k(V,\Omega_I)_\tau^\tau \le c\|V\|_\tau^\tau
\le c|J|\|g\|_{\BMO}^\tau + c|J|\big(|J|\|U'\|_{L^\infty(J^\star)}\big)^\tau
\end{equation}
as in \eqref{est-V} with the interval $3J$ replaced with $[x_{j_1}-|J|/2,x_{j_1}]\cup J\subset J^\star$.

For the estimation of $U'$ we use \eqref{est-Up}, which together with \eqref{est-V5} gives \eqref{est_E_I}.
Estimate \eqref{est_E_I} is also valid in the case when $J_{j_1}$ is unbounded because $g=U=0$ here.
Thus, \eqref{est_E_I} is proved.

Using estimates \eqref{est_E_I} we obtain
\begin{multline}\label{est_E_Im}
\sum_{I\in\II_m}|I|^{-1}E_k(g,\Omega_I)_\tau^\tau\\
\le c\sum_{I\in\II_m}|I|^{-1}\bigg(\sum_{J_\nu\subset\Omega_I} |J_\nu|
+\sum_{\substack{0<|J_\mu\cap\Omega_I|<|\Omega_I|\\J_\mu\setminus\Omega_I\ne\emptyset}}
\min\Big\{|J_\mu|,\frac{t_m^{1+\tau}}{|J_\mu|^\tau}\Big\}\bigg)\|g\|_{\BMO}^\tau\\
\le c \bigg(\sum_{|J_\nu|\le t_m}\frac{|J_\nu|}{t_m} + \sum_{|J_\nu|> t_m}\frac{t_m^{\tau}}{|J_\mu|^\tau}\bigg)\|g\|_{\BMO}^\tau.
\end{multline}
In the last inequality we use that every $J_\nu$ with $|J_\nu|\le t_m$ may belong to at most $2k-1$ different intervals $\Omega_I$, $I\in\II_m$,
and that every $J_\mu$ (independently of $|J_\nu|\le t_m$ or $t_m<|J_\nu|<\infty $) may partially cover at most $4k-2$ different intervals $\Omega_I$, $I\in\II_m$.

Finally, taking a sum on $m$ in \eqref{est_E_Im} and using $t_{m+1}\le \rho t_m$ and $N\le (k+1)n$ we obtain
\begin{multline}\label{est_E_Imz}
\sum_{I\in\II}|I|^{-1}E_k(g,\Omega_I)_\tau^\tau
\le c \sum_{m\in\ZZ}\bigg(\sum_{|J_\nu|\le t_m}\frac{|J_\nu|}{t_m} + \sum_{|J_\nu|> t_m}\frac{t_m^{\tau}}{|J_\mu|^\tau}\bigg)\|g\|_{\BMO}^\tau\\
= c \sum_{\nu=1}^N\bigg(\sum_{\substack{m\in\ZZ\\ |J_\nu|\le t_m}}\frac{|J_\nu|}{t_m}
+ \sum_{\substack{m\in\ZZ\\ |J_\nu|> t_m}}\frac{t_m^{\tau}}{|J_\mu|^\tau}\bigg)\|g\|_{\BMO}^\tau
\le c n \|g\|_{\BMO}^\tau.
\end{multline}
In view of \eqref{B-tau} this completes the proof of the theorem in the case $1\le\tau<\infty$.

(b) Let $0<\tau<1$. We shall use the identification $\BB_\tau^{\alpha,k}=\BB_\tau^{\alpha,k}(E,1)$, see \eqref{def-B-q}.
From \eqref{est_E_I} with $\tau=1$ and the concavity of $y^\tau$ for $0<\tau<1$ we obtain for every $I\in\II_m$, $m\in\ZZ$, the inequality
\begin{multline}\label{est_E_I1}
(|I|^{-1}E_k(g,\Omega_I)_1)^\tau\\
\le c \bigg(\sum_{J_\nu\subset\Omega_I}\frac{|J_\nu|^\tau}{t_m^\tau}
+\sum_{\substack{0<|J_\mu\cap\Omega_I|<|\Omega_I|\\J_\mu\setminus\Omega_I\ne\emptyset}}
\min\Big\{\frac{|J_\mu|^\tau}{t_m^\tau},\frac{t_m^{\tau}}{|J_\mu|^\tau}\Big\}\bigg)\|g\|_{\BMO}^\tau.
\end{multline}
Now, proceeding as in the proof of \eqref{est_E_Im} and \eqref{est_E_Imz} we obtain from \eqref{est_E_I1}
\begin{multline*}
\sum_{I\in\II}(|I|^{-1}E_k(g,\Omega_I)_1)^\tau\\
\le c \sum_{\nu=1}^N\bigg(\sum_{\substack{m\in\ZZ\\ |J_\nu|\le t_m}}\frac{|J_\nu|^\tau}{t_m^\tau}
+ \sum_{\substack{m\in\ZZ\\ |J_\nu|> t_m}}\frac{t_m^{\tau}}{|J_\mu|^\tau}\bigg)\|g\|_{\BMO}^\tau
\le c n \|g\|_{\BMO}^\tau.
\end{multline*}
In view of \eqref{def-B-q} this completes the proof of the theorem in the case $0<\tau<1$.
$\hfill \square$

\section{Appendix}\label{sec:appendix}

\subsection{Proofs of Lemma~\ref{lem:embed-BMO} and Theorem~\ref{thm:embed-B-BMO}}\label{subsec:appendix1}

\noindent
{\it Proof of Lemma~\ref{lem:embed-BMO}.}
Let $\{\cc_\Q\}_{\Q\in\QQ}$ be a sequence of complex numbers
and fix an compact interval $J\subset \RR$.
Consider first the case when $0<\tau<p$. By \cite[Theorem~3.3]{KP} we have
\begin{equation}\label{embed-1}
\Big\|\sum_{\Q\in\QQ, \Q\subset J}|\cc_\Q\varphi_\Q|\Big\|_p
\le c\Big(\sum_{\Q\in\QQ, \Q\subset J}\|\cc_\Q\varphi_\Q\|^\tau_p\Big)^{1/\tau}.
\end{equation}
Clearly,
$\|\cc_\Q\varphi_\Q\|_p
\le  c|\cc_\Q| |\Q|^{1/p}
\le c|\cc_\Q| |J|^{1/p},
$
which along with \eqref{embed-1} implies \eqref{embed-BMO}.

In the case $\tau \ge p$ we choose $q>\tau$ and use H\"{o}lder's inequality
and \eqref{embed-BMO} in the proven case from above to obtain
\begin{align*}
\Big\|\frac{1}{|J|^{1/p}}\sum_{\substack{\Q\in\QQ,\\ \Q\subset J}}|\cc_\Q\varphi_\Q|\Big\|_p
\le \Big\|\frac{1}{|J|^{1/q}}\sum_{\substack{\Q\in\QQ,\\ \Q\subset J}}|\cc_\Q\varphi_\Q|\Big\|_q
\le c\Big(\sum_{\substack{\Q\in\QQ,\\ \Q\subset J}} |\cc_\Q|^\tau\Big)^{1/\tau}.
\end{align*}
The proof is complete.
$\hfill \square$

\smallskip
\noindent
{\it Proof of Theorem~\ref{thm:embed-B-BMO}.}
Part (b) is trivial. For the proof of part (a) assume $\tau>1$.
Denote
\begin{equation*}
f_\nu:=\sum_{j>\nu}\sum_{\Q\in\QQ_j} \cc_\Q\varphi_\Q,
\quad \nu\in\ZZ.
\end{equation*}
We claim that
\begin{equation}\label{est-BMO-fk}
\|f_\nu\|_{\BMO} \le c\Big(\sum_{j>\nu}\sum_{\Q\in\QQ_j} |\cc_\Q|^\tau\Big)^{1/\tau}
=:\|\{\cc_\Q\}\|_{\ell^\tau(\nu)}.
\end{equation}
Let $J$ be an arbitrary compact interval in $\RR$.
Then there exist $m\in\ZZ$ such that
if $\Q\in\QQ_m$ and $\Q\cap J\ne \emptyset$, then $|\Q|\sim |J|$  and $\Q\subset 2J$.
Denote $\tilde{J}:=2J$.

First, we consider the less favorable case $\nu<m$.
We split $f_\nu$ into two: $f_\nu=f_m+(f_\nu-f_m)$.
Using Lemma~\ref{lem:embed-BMO} we get
\begin{multline}\label{BMO-fm}
\frac{1}{|J|}\int_J|f_m(x)|dx
\le \frac{c}{|\tilde{J}|}\int_{\tilde{J}}\sum_{j>m}\sum_{\substack{\Q\in\QQ_j,\\\Q\subset \tilde{J}}} |\cc_\Q\varphi_\Q(x)|dx\\
\le c\Big(\sum_{j>m}\sum_{\Q\in\QQ_j} |\cc_\Q|^\tau\Big)^{1/\tau}
\le c\|\{\cc_\Q\}\|_{\ell^\tau(\nu)}.
\end{multline}

Denote $F_{\nu m}:=f_\nu-f_m$ and fix $y\in J$. We claim that
\begin{equation}\label{est-Fkm}
\frac{1}{|J|}\int_J|F_{\nu m}(x)-F_{\nu m}(y)|dx
\le c\|\{\cc_\Q\}\|_{\ell^\tau(\nu)}.
\end{equation}
Indeed, let $\Q\in\QQ_j$, $j\le m$,
and assume $\Q\cap J\ne \emptyset$.
Then for $x\in \Q$
\begin{equation*}
|\varphi_\Q(x)- \varphi_\Q(y)| \le |x-y|\|\varphi_\Q'\|_\infty \le c|x-y||\Q|^{-1}.
\end{equation*}
Fix $x\in J$ and assume that $x$ belongs to the interior of some $\Q^\star\in \QQ_m$.
Using the above we get
\begin{align*}
|F_{\nu m}(x)-F_{\nu m}(y)|
&\le \sum_{j=\nu +1}^m\sum_{\substack{\Q\in\QQ_j,\\ \Q\ni x}} |\cc_\Q||\varphi_\Q(x)-\varphi(y)|
\\
&\le c\|\{\cc_\Q\}\|_{\ell^\tau(\nu)} \sum_{j=\nu +1}^m\sum_{\substack{\Q\in\QQ_j,\\ \Q\ni x}} |x-y||\Q|^{-1}
\\
&\le c\|\{\cc_\Q\}\|_{\ell^\tau(\nu)} |J|\sum_{j\le m}\sum_{\substack{\Q\in\QQ_j,\\ \Q\ni x}}|\Q|^{-1}
\\
&\le c\|\{\cc_\Q\}\|_{\ell^\tau(\nu)} |J||\Q^\star|^{-1}
\le c\|\{\cc_\Q\}\|_{\ell^\tau(\nu)}.
\end{align*}
Here we used that
$
\sum_{j\le m}\sum_{\Q\in\QQ_j, \Q\ni x}|\Q|^{-1} \le c|\Q^\star|^{-1},
$
which follows from the conditions on the underlying regular multilevel partition $\II$.
Estimate \eqref{est-Fkm} follows readily from the above inequalities.

From \eqref{BMO-fm} and \eqref{est-Fkm} it follows that
\begin{equation*}
\frac{1}{|J|}\int_J|f_\nu(x)-F_{\nu m}(y)|dx
\le c\|\{\cc_\Q\}\|_{\ell^\tau(\nu)},
\end{equation*}
which implies \eqref{est-BMO-fk}.

In the easier case $\nu\ge m$ \eqref{est-BMO-fk} will follow directly from an estimate similar to \eqref{BMO-fm}.
In turn, \eqref{est-BMO-fk} implies that for any $\nu, \mu\in\ZZ$, $\mu>\nu$,
\begin{equation*}
\|f_\nu-f_\mu\|_{\BMO} \le c\Big(\sum_{j=\nu+1}^\mu\sum_{\Q\in\QQ_j} |\cc_\Q|^\tau\Big)^{1/\tau} \to 0
\quad\hbox{as}\quad \nu, \mu\to -\infty.
\end{equation*}
Since $\BMO$ is complete, it follows that $\lim_{\nu\to-\infty} f_\nu =f$ for some $f\in \BMO$,
where the convergence is in the $\BMO$-norm.
It also follows that
$\|f\|_{\BMO} \le c\|\{\cc_\Q\}\|_{\ell^\tau}$,
which confirms \eqref{embed-B-BMO}.

Finally, because the norm $\|\{\cc_\Q\}\|_{\ell^\tau}$ does not change when reshuffling the terms in its definition,
it readily follows from the above proof that the convergence
in $\sum_{\Q\in\QQ} \cc_\Q\varphi_\Q$ is unconditional in $\BMO$.
$\hfill \square$

\subsection{Proof of equivalence \eqref{B-tau}}\label{subsec:appendix2}

Denote
\begin{equation*}
\|f\|_{\BB^{\alpha, k}_\tau(E)} := \Big(\sum_{I \in \II} \big(|I|^{-\alpha}\omega_k(f, \Om_I)_\tau)^\tau\Big)^{1/\tau}.
\end{equation*}
Denote by $\sD_m$ the $m$th level dyadic intervals
($|J|=2^{-m}$ if $J\in\sD_m$)
and set $\sD:=\cup_{m\in\ZZ}\sD_m$.
Clearly, see \eqref{omega-J},
\begin{equation*}
\omega_k(f, 2^{-m})_\tau^\tau \le \sum_{J\in\sD_m} \omega_k(f, (2k+1)J)_\tau^\tau.
\end{equation*}
From the conditions on $\II$ it follows that for each $J\in\sD$
the interval $(2k+1)J$ is contained in some interval $\Omega_I$, $I\in\II$, of minimum lenght (hence, $|\Omega_I|\sim |J|$),
and each $\Omega_I$, $I\in\II$, contains a uniformly bounded number of such intervals $J\in\sD$.
Therefore,
\begin{multline*}
\|f\|_{\BB^{\alpha, k}_\tau}^\tau \sim \sum_{m\in\ZZ}2^{m}\omega_k(f, 2^{-m})_\tau^\tau\\
\le c\sum_{J\in\sD} |J|^{-1}\omega_k(f, (2k+1)J)_\tau^\tau
\le c \sum_{I \in \II} \big(|I|^{- \frac{1}{\tau}}\omega_k(f, \Om_I)_\tau)^\tau
\end{multline*}
and hence
$\|f\|_{\BB^{\alpha, k}_\tau} \le c \|f\|_{\BB^{\alpha, k}_\tau(E)}$.

For the estimate in the other direction we use \eqref{omega-aver}.
We obtain
\begin{align*}
\sum_{I \in \II, \frac{1}{2}<2^m|I|\le 2} \big(|I|^{- \frac{1}{\tau}}&\omega_k(f, \Om_I)_\tau)^\tau\\
&\le c\sum_{I \in \II, \frac{1}{2}<2^m|I|\le 2} |I|^{-2}\int_0^{|\Omega_I|}\int_{\Omega_I}|\Delta^k_h(f,x,\Omega_I)|^\tau dxdh
\\
& \le c2^{2m} \int_0^{c2^{-m}}\int_\RR |\Delta^k_hf (x)|^\tau dxdh
\\
&\le c 2^m\omega_k(f, c2^{-m})_\tau^\tau\le c 2^m\omega_k(f, 2^{-m})_\tau^\tau.
\end{align*}
Here we used that only finitely many of the intervals $\{\Omega_I: I \in \II, \frac{1}{2}\!<\!2^m|I|\!\le\!2\}$
may overlap at any point $x\in\RR$,
and $\omega_k(f, c2^{-m})_\tau \le c'\omega_k(f, 2^{-m})_\tau$.
From above and \eqref{discr-B-mogulus} we get
\begin{align*}
\|f\|_{\BB^{\alpha, k}_\tau(E)}=\sum_{I \in \II} \big(|I|^{- \frac{1}{\tau}}\omega_k(f, \Om_I)_\tau)^\tau
\le c\sum_{m\in\ZZ} 2^m\omega_k(f, 2^{-m})_\tau^\tau
\le c\|f\|^\tau_{\BB^{\alpha, k}_\tau}.
\end{align*}
This and the estimate in the other direction from above yield the equivalence
$\|f\|_{\BB^{\alpha, k}_\tau} \sim \|f\|_{\BB^{\alpha, k}_\tau(E)}$.
$\hfill \square$

\subsection{Proof of Theorem~\ref{thm:rep-B}}\label{subsec:appendix3}

Let $f\in \BB^{\alpha, k}_\tau(E, q)$.
In light of Proposition~\ref{prop:BMO-embed} there exists a polynomial $P\in\Pi_k$ such that
$\|f-P\|_{\BMO} \le c \|f\|_{\BB^{\alpha, k}_\tau}(E, q)$.
Let $T_{m, q}$ be the quasi-interpolant from \eqref{def-proj-2}

(a) First we show that
\begin{equation}\label{Q-tau-to-f}
\lim_{m\to \infty} \|f-P-T_{m,q}(f-P)\|_{\BMO} = 0.
\end{equation}
Fix $\eps>0$. In light of \eqref{def-B-q} there exists $m_0\in \NN$ such that
\begin{equation}\label{B-eps}
\sum_{j=m_0}^\infty \sum_{I\in\II_j} |I|^{-\tau/q}E_k(f, \Omega_I)_q^\tau <\eps^{\tau}.
\end{equation}
Fix $m\ge m_0$.
Let $J$ be an arbitrary compact interval and let $\nu$ be its level (see \S\ref{subsec:nested_partitions}).
We next consider two cases depending on the size of  $|J|$.

\smallskip

{\em Case 1:} $\nu> m$.
There exist two adjacent intervals $I_1$, $I_2$ in $\II_\nu$ such that
$J\subset I_1\cup I_2$, $|J|\sim |I_1|\sim |I_2|$.
For an appropriate constant $c^\diamond$ (to be selected) we have
\begin{align}\label{Q-to-f-1}
\frac{1}{|J|}&\int_J \big|f(x)-P(x)- T_{m,q}(f-P)(x)-c^\diamond\big|^q dx
\\
&= \frac{1}{|J|}\int_J |f(x)-T_{m,q}(f)(x)-c^\diamond|^q dx
 \le \frac{c}{|J|}\int_J |f(x)-T_{\nu,q}(f)(x)|^q dx \nonumber
\\
&+ c\big\|T_{\nu,q}(f)-T_{m,q}(f)-c^\diamond\big\|_{L^\infty(J)}^q
=: S_1+S_2. \nonumber
\end{align}
To estimate $S_1$ we use Lemma~\ref{lem:quasi-int} and \eqref{B-eps} to obtain
\begin{multline}\label{est-S1}
S_1\le \frac{1}{|J|}\int_{I_1\cup I_2}|f(x)-T_{\nu, q}(f)(x)|^q dx\\
\le c|I_1|^{-1}E_k(f, \Omega_{I_1})_q^q + c|I_2|^{-1}E_k(f, \Omega_{I_2})_q^q <c\eps^{q}.  
\end{multline}
To estimate $S_2$ we shall use the abbreviated notation $\qq_j:=T_{j,q}(f)-T_{j-1,q}(f)$ (see \eqref{def:qt}).
We fix $y\in J$ and select the constant
$c^\diamond:= T_{\nu,q}(f)(y)- T_{m,q}(f)(y)$.
Then for any $x\in J$ we have
\begin{align*}
|T_{\nu, q}(f)(x)-T_{m,q}(f)(x)-c^\diamond|
=\Big|\!\!\sum_{j=m+1}^\nu (\qq_j(x)-\qq_j(y))\Big|
\le |J|\!\!\sum_{j=m+1}^\nu\|\qq_j'\|_{L^\infty(J)}.
\end{align*}
The choice of $\nu$ implies that for any $j=m+1, \dots, \nu$
there exist  two adjacent intervals $I_j',I_j''$ in $\II_j$ such that $J\subset I_j'\cup I_j''$.
Using that $\qq_j$ is a polynomial of degree $k-1$ on $I_j'$ and on $I_j''$ we obtain from \lemref{lem:poly-norms} and \eqref{cond-rho}
\begin{align*}
|J|\sum_{j=m+1}^\nu\|\qq_j'\|_{L^\infty(J)} &\le |J|\sum_{j=m+1}^\nu\Big(\|\qq_j'\|_{L^\infty(I_j')}+\|\qq_j'\|_{L^\infty(I_j'')}\Big)
\\
&\le c|J|\sum_{j=m+1}^\nu(|I_j'|^{-1}\|\qq_j\|_{L^\infty(I_j')} + |I_j''|^{-1}\|\qq_j\|_{L^\infty(I_j'')})
\\
&\le c\sum_{j=m+1}^\nu\frac{|J|}{|I_j'|}
\big(|I_j'|^{-1/q}\|\qq_j\|_{L^q(I_j')} + |I_j''|^{-1/q}\|\qq_j\|_{L^q(I_j'')}\big).
\end{align*}
Using \eqref{est-qj}, \eqref{B-eps} and \eqref{r-rho} in the above, we obtain
\begin{align*}
S_2=c\big\|T_{\nu, q}(f)-T_{m,q}(f)-c^\diamond\big\|_{L^\infty(J)}^q <c\eps^q.
\end{align*}
This together with \eqref{Q-to-f-1} and \eqref{est-S1} implies
\begin{equation}\label{f-P-Q}
\frac{1}{|J|}\int_J \big|f(x)-P(x)- T_{m,q}(f-P)(x)-c^\diamond\big|^q dx \le c\eps^{q}.
\end{equation}

\smallskip

{\em Case 2:} $\nu\le m$. Hence $|J| \ge c|I|$ for all $I\in\II_m$ and $\sum_{I\in\II_m, I\cap J\ne \emptyset}|I|\le c|J|$.
Using \lemref{lem:quasi-int} and \eqref{B-eps} we obtain
\begin{multline*}
\frac{1}{|J|}\int_J |f(x)-T_{m,q}(f)(x)|^q dx
\le \frac{1}{|J|}\sum_{I\in\II_m, I\cap J\ne \emptyset}\int_I|f(x)-T_{m,q}(f)(x)|^q dx \nonumber
\\
\le \frac{c}{|J|}\sum_{I\in\II_m, I\cap J\ne \emptyset} E_k(f, \Omega_I)_q^q
\le \frac{c}{|J|}\sum_{I\in\II_m, I\cap J\ne \emptyset}|I| \eps^q \le c \eps^q.
\end{multline*}
In turn, this and estimate \eqref{f-P-Q} yield
\begin{equation*}
\|f-P-T_{m,q}(f-P)\|_{\BMO} \le c\eps, \quad \forall m\ge m_0,
\end{equation*}
which implies \eqref{Q-tau-to-f}.

\smallskip

(b) We next prove that
\begin{equation}\label{conv-Q-m-tau}
\lim_{m\to-\infty}\|T_{m,q}(f-P)\|_{\BMO}=0.
\end{equation}
Let $\eps>0$. By \eqref{def-B-q} it follows that there exists $m_1\in \ZZ$ such that
\begin{equation}\label{B-eps-2}
\sum_{j=-\infty}^{m_1} \sum_{I\in\II_j} |I|^{-\tau/q}E_k(f, \Omega_I)_q^\tau <\eps^{\tau}.
\end{equation}
Fix $m< m_1$.
Let $J$ be an arbitrary compact interval and let $\nu-1$ be its level (see \S\ref{subsec:nested_partitions}).
Then $J$ contains some interval $I\in \II_\nu$ and $|J|\sim |I|$.

As in part (a) we shall use the abbreviated notation $\qq_j:=T_{j,q}(f)-T_{j-1,q}(f)$.
Observe that $T_{m,q}(f-P)=T_{m,q}(f)-P$.
Using this we write
\begin{align*}
T_{m,q}(f-P) = \!\!\sum_{j=N+1}^m \qq_j + T_{N,q}(f)-P
= \!\!\sum_{j=N+1}^m\sum_{\Q\in\QQ_j}b_{Q,q}(f)\varphi_\Q + T_{N,q}(f)-P,
\end{align*}
where $N<m$, $N<\nu$ and $\nu-N$ is sufficiently large (to be determined).
Clearly, for any constant $c^\star$ (to be selected) there exists a constant $c^{\star\star}$ such that
\begin{align}\label{Qm-1}
\frac{1}{|J|}\int_J|T_{m,q}(f-P)-c^{\star\star}|^q dx
&\le c\Big\|\sum_{j=N+1}^m\sum_{\Q\in\QQ_j}b_{Q,q}(f)\varphi_\Q\Big\|_{\BMO}^q 
\\
&+ \frac{c}{|J|}\int_J |T_{N,q}(f)(x)-P(x)-c^\star|^q dx
=: S_1+S_2.\nonumber
\end{align}
To estimate $S_1$ we invoke Theorem~\ref{thm:embed-B-BMO}, \eqref{NQ=b-norm1}, \eqref{NQ=b-normt}, \eqref{est-qj}, \eqref{B-eps-2} and obtain
\begin{align}\label{S-1}
S_1&\le c\Big(\sum_{j=-\infty}^m \sum_{\Q\in\QQ_j} |b_{Q,q}(f)|^\tau\Big)^{q/\tau}
\le c\Big(\sum_{j=-\infty}^m \sum_{I\in\II_j}|I|^{-\tau/q}\|\qq_j\|_{L^q(I)}^\tau\Big)^{q/\tau}
\\
&\le c\Big(\sum_{j=-\infty}^{m} \sum_{I\in\II_j} |I|^{-\tau/q}E_k(f, \Omega_I)_q^\tau\Big)^{q/\tau} \nonumber
 <c\eps^{q}. \nonumber
\end{align}

To estimate $S_2$ we recall that $N<\nu$ and hence
there are two adjacent intervals $I_1, I_2$ in $\II_N$ such that $J\subset I_1\cup I_2$.
Let $I^\diamond\in \II_{N-1}$ be the only parent of $I_1$ ($I_1\subset I^\diamond$).
Clearly, $\Omega_{I_1}\cup \Omega_{I_2} \subset \Omega_{I^\diamond}$.
Let $R\in\Pi_k$ be a polynomial such that
\begin{equation}\label{def-R}
\|f-R\|_{L^q(\Omega_{I^\diamond})} \le cE_k(f, \Omega_{I^\diamond})_q.
\end{equation}
We now choose the constant $c^\star$ to be $c^\star:= R(y)-P(y)$, where $y\in J$ is fixed.
We have
\begin{align}\label{S-2}
 &S_2 \le \frac{c}{|J|}\int_J|T_{N,q}(f)(x)-P(x)-c^\star|^q dx
\\
&\le \frac{c}{|J|}\int_J |T_{N}(\sP_{N,q}-R)(x)|^q dx + \frac{c}{|J|}\int_J |R(x)-P(x)-c^\star|^q dx
=: U_1+U_2. \nonumber
\end{align}
Using \eqref{Q-m-1} and \eqref{equiv-norms} we get
\begin{align*}
U_1&= \frac{c}{|J|}\int_J |T_{N}(\sP_{N,q}-R)(x)|^q dx
\le c\|T_{N}(\sP_{N,q}-R)\|_{L^\infty(I_1\cup I_2)}^q
\\
&\le c\|\sP_{N,q}-R\|_{L^\infty(\Omega_{I_1}\cup \Omega_{I_2})}^q
\le c\max_{I\in\II_N, I\subset \Omega_{I_1}\cup \Omega_{I_2}}\|\sP_{N,q}-R\|_{L^\infty(I)}^q
\displaybreak \\
&\le c\sum_{I\in\II_N, I\subset \Omega_{I_1}\cup \Omega_{I_2}}|I|^{-1}\|\sP_{N,q}-R\|_{L^q(I)}^q
\\
&\le c\sum_{I\in\II_N, I\subset \Omega_{I_1}\cup \Omega_{I_2}}
\big(|I|^{-1}\|f-\sP_{N,q}\|_{L^q(I)}^q + |I|^{-1}\|f-R\|_{L^q(I)}^q\big).
\end{align*}
By the definition of $\sP_{N,q}$ we have
$\|f-\sP_{N,q}\|_{L^q(I)} \le cE_k(f, I)_q \le cE_k(f, \Omega_I)_q$
(see \eqref{def-proj} and \eqref{project}).
We use this, \eqref{def-R}, \eqref{B-eps-2} and the above estimates for $U_1$ to obtain
\begin{equation}\label{U-1}
U_1\le c\sum_{j=N-1}^N \sum_{I\in\II_j, I\cap (\Omega_{I_1}\cup\Omega_{I_2})\ne \emptyset} |I|^{-1}E_k(f, \Omega_I)_q^q <c\eps^{q}.
\end{equation}
Here we used that in the double sum above there is a constant number (depending only on $k$) of terms.

In estimating $U_2$ we shall use the abbreviated notation $I_\star:=I_1\cup I_2$, where $I_1,I_2\in\II_N$ are determined above.
Using that $R$ and $P$ are polynomials on $I_\star$ and \lemref{lem:poly-norms} we get
\begin{align*}
U_2&=\frac{c}{|J|}\int_J |R(x)-P(x)-(R(y)-P(y))|^q dx
\le c(|J|\|(R-P)'\|_{L^\infty(I_\star)})^q
\\
& = c\big(|J|\|(R-P -\tilde{c})'\|_{L^\infty(I_\star)}\big)^q
\le c(|J|/|I_\star|)^q \|R-P -\tilde{c}\|_{L^\infty(I_\star)}^q
\\
&\le c(|J|/|I_\star|)^q |I_\star|^{-1}\|R-P -\tilde{c}\|_{L^q(I_\star)}^q,
\end{align*}
where the constant $\tilde{c}$ is defined by $\tilde{c}:= \avg_{I_\star}(f-P)$.
We now use \eqref{def-R}, \eqref{equiv-norms-BMO}, \eqref{B-eps-2}, and obtain
\begin{align*}
|I_\star|^{-1}&\|R-P -\tilde{c}\|_{L^q(I_\star)}^q
\\
&\le c|I_\star|^{-1}\|f-R\|_{L^q(I_\star)}^q
+ c|I_\star|^{-1}\|f-P -\avg_{I_\star}(f-P)\|_{L^q(I_\star)}^q
\\
&\le c|I^\diamond|^{-1}E_k(f, \Omega_{I^\diamond})_q^q
+ c\|f-P\|_{\BMO}^q
\le c\eps^{q} + c\|f-P\|_{\BMO}^q.
\end{align*}
On the other hand, because $|J|\sim|I|$ with $I\in \II_\nu$ and
$I_\star=I_1\cup I_2$ with $I_1, I_2\in \II_N$, we infer from \eqref{r-rho}
that
$|J|/|I_\star| \le c\rho^{\nu-N}$.
Putting all of the above together we obtain
\begin{equation*}
U_2 \le c\rho^{(\nu-N)q}(\eps^{q} + \|f-P\|_{\BMO}^q).
\end{equation*}
Combining this with \eqref{S-2} and \eqref{U-1} we get
\begin{equation*}
S_2 \le c\eps^{q} + c\rho^{(\nu-N)q}\|f-P\|_{\BMO}^q.
\end{equation*}
In turn, this along with \eqref{Qm-1} and \eqref{S-1} yield
\begin{equation*}
\frac{1}{|J|}\int_J|T_{m,q}(f-P)-c^\star|^q dx \le c\eps^{q} + c\rho^{(\nu-N)q}\|f-P\|_{\BMO}^q,
\quad\forall m<m_1.
\end{equation*}
Since the constant $c$ in this estimate is independent of $N$ and $f-P\in\BMO$,
by letting $N\to -\infty$ we arrive at
\begin{equation*}
\frac{1}{|J|}\int_J|T_{m,q}(f-P)-c^\star|^q dx \le c\eps^{q},
\quad\forall m<m_1.
\end{equation*}
This estimate implies
$\|T_{m,q}(f-P)\|_{\BMO} \le c\eps$ for all $m<m_1$,
which yields \eqref{conv-Q-m-tau}.

Clearly, decomposition \eqref{rep-B} follows at once by \eqref{Q-tau-to-f} and \eqref{conv-Q-m-tau}.
Inequality \eqref{BE-embed-BMO} follows by Lemma~\ref{lem:BQ-BE}.
The unconditional convergence in \eqref{rep-B} is a consequence of Theorem~\ref{thm:embed-B-BMO}.
Finally, the unconditional convergence in $\BMO$ of the series in \eqref{rep-B}
and the fact that each $\varphi_Q$ is in $C_0(\RR)$
leads to the conclusion that $f-P$ is in $\VMO$.
$\hfill \square$

\end{sloppypar}

\end{document}